\documentclass[12pt]{amsart}

\usepackage[centertags]{amsmath}
\usepackage{amsthm,amsfonts}

\usepackage{amssymb}
\usepackage{epsfig}
\usepackage{amssymb,latexsym}
\usepackage{indentfirst}

\begin{document}
\theoremstyle{plain}
\newtheorem{theorem}{Theorem}[section]
\newtheorem{lemma}[theorem]{Lemma}
\newtheorem{corollary}[theorem]{Corollary}
\newtheorem{proposition}[theorem]{Proposition}
\newtheorem{example}[theorem]{Example}
\newtheorem{examples}[theorem]{Examples}
\theoremstyle{definition}
\newtheorem{notations}[theorem]{Notations}
\newtheorem{notation}[theorem]{Notation}
\newtheorem{remark}[theorem]{Remark}
\newtheorem{remarks}[theorem]{Remarks}
\newtheorem{definition}[theorem]{Definition}
\newtheorem{claim}[theorem]{Claim}
\newtheorem{assumption}[theorem]{Assumption}
\numberwithin{equation}{section}
\newtheorem{examplerm}[theorem]{Example}
\newtheorem{propositionrm}[theorem]{Proposition}


\newcommand{\binomial}[2]{\left(\begin{array}{c}#1\\#2\end{array}\right)}
\newcommand{\zar}{{\rm zar}}
\newcommand{\an}{{\rm an}}
\newcommand{\red}{{\rm red}}
\newcommand{\codim}{{\rm codim}}
\newcommand{\rank}{{\rm rank}}
\newcommand{\Pic}{{\rm Pic}}
\newcommand{\Div}{{\rm Div}}
\newcommand{\Hom}{{\rm Hom}}
\newcommand{\Ima}{{\rm Im}}
\newcommand{\Ker}{{\rm Ker}}
\newcommand{\Spec}{{\rm Spec}}
\newcommand{\sing}{{\rm sing}}
\newcommand{\reg}{{\rm reg}}
\newcommand{\Char}{{\rm char}}
\newcommand{\Tr}{{\rm Tr}}
\newcommand{\tr}{{\rm tr}}
\newcommand{\supp}{{\rm supp}}
\newcommand{\Gal}{{\rm Gal}}
\newcommand{\Min}{{\rm Min \ }}
\newcommand{\Max}{{\rm Max \ }}
\newcommand{\soplus}[1]{\stackrel{#1}{\oplus}}
\newcommand{\dlog}{{\rm dlog}\,}    
\newcommand{\limdir}[1]{{\displaystyle{\mathop{\rm
lim}_{\buildrel\longrightarrow\over{#1}}}}\,}
\newcommand{\liminv}[1]{{\displaystyle{\mathop{\rm
lim}_{\buildrel\longleftarrow\over{#1}}}}\,}
\newcommand{\boxtensor}{{\Box\kern-9.03pt\raise1.42pt\hbox{$\times$}}}
\newcommand{\sext}{\mbox{${\mathcal E}xt\,$}}
\newcommand{\shom}{\mbox{${\mathcal H}om\,$}}
\newcommand{\coker}{{\rm coker}\,}
\renewcommand{\iff}{\mbox{ $\Longleftrightarrow$ }}
\newcommand{\onto}{\mbox{$\,\>>>\hspace{-.5cm}\to\hspace{.15cm}$}}

\newenvironment{pf}{\noindent\textbf{Proof.}\quad}{\hfill{$\Box$}}

\newcommand{\sA}{{\mathcal A}}
\newcommand{\sB}{{\mathcal B}}
\newcommand{\sC}{{\mathcal C}}
\newcommand{\sD}{{\mathcal D}}
\newcommand{\sE}{{\mathcal E}}
\newcommand{\sF}{{\mathcal F}}
\newcommand{\sG}{{\mathcal G}}
\newcommand{\sH}{{\mathcal H}}
\newcommand{\sI}{{\mathcal I}}
\newcommand{\sJ}{{\mathcal J}}
\newcommand{\sK}{{\mathcal K}}
\newcommand{\sL}{{\mathcal L}}
\newcommand{\sM}{{\mathcal M}}
\newcommand{\sN}{{\mathcal N}}
\newcommand{\sO}{{\mathcal O}}
\newcommand{\sP}{{\mathcal P}}
\newcommand{\sQ}{{\mathcal Q}}
\newcommand{\sR}{{\mathcal R}}
\newcommand{\sS}{{\mathcal S}}
\newcommand{\sT}{{\mathcal T}}
\newcommand{\sU}{{\mathcal U}}
\newcommand{\sV}{{\mathcal V}}
\newcommand{\sW}{{\mathcal W}}
\newcommand{\sX}{{\mathcal X}}
\newcommand{\sY}{{\mathcal Y}}
\newcommand{\sZ}{{\mathcal Z}}

\newcommand{\A}{{\mathbb A}}
\newcommand{\B}{{\mathbb B}}
\newcommand{\C}{{\mathbb C}}
\newcommand{\D}{{\mathbb D}}
\newcommand{\E}{{\mathbb E}}
\newcommand{\F}{{\mathbb F}}
\newcommand{\G}{{\mathbb G}}
\newcommand{\HH}{{\mathbb H}}
\newcommand{\I}{{\mathbb I}}
\newcommand{\J}{{\mathbb J}}
\newcommand{\M}{{\mathbb M}}
\newcommand{\N}{{\mathbb N}}
\renewcommand{\P}{{\mathbb P}}
\newcommand{\Q}{{\mathbb Q}}
\newcommand{\R}{{\mathbb R}}
\newcommand{\T}{{\mathbb T}}
\newcommand{\U}{{\mathbb U}}
\newcommand{\V}{{\mathbb V}}
\newcommand{\W}{{\mathbb W}}
\newcommand{\X}{{\mathbb X}}
\newcommand{\Y}{{\mathbb Y}}
\newcommand{\Z}{{\mathbb Z}}

\newcommand{\be}{\begin{eqnarray}}
\newcommand{\ee}{\end{eqnarray}}
\newcommand{\nn}{{\nonumber}}
\newcommand{\dd}{\displaystyle}
\newcommand{\ra}{\rightarrow}

\title[Improved asymptotic bounds for codes]{Improved asymptotic bounds for codes using distinguished divisors of global function fields}

\author{Harald Niederreiter and Ferruh \"Ozbudak}

\maketitle

\begin{center}
Harald Niederreiter \\
Department of Mathematics, National University of Singapore \\
2 Science Drive 2, Singapore 117543, Republic of Singapore \\
e-mail: nied@math.nus.edu.sg
\\$\mbox{}$ \\
Ferruh \"{O}zbudak\\
Temasek Laboratories,
National University of Singapore\\
5 Sports Drive 2, 117508 Singapore, Republic of Singapore \\
and \\
Department of Mathematics, Middle East Technical University \\
\.{I}n\"on\"u Bulvar{\i}, 06531, Ankara, Turkey \\
e-mail: ozbudak@metu.edu.tr
\end{center}

\abstract
For a prime power $q$, let $\alpha_q$ be the standard function in the
asymptotic theory of codes, that is, $\alpha_q(\delta)$ is the largest
asymptotic information rate that can be achieved for a given asymptotic
relative minimum distance $\delta$ of $q$-ary codes. In recent years the
Tsfasman-Vl\u{a}du\c{t}-Zink lower bound on $\alpha_q(\delta)$ was
improved by Elkies, Xing, and Niederreiter and \"Ozbudak. In this paper
we show further improvements on these bounds by using distinguished
divisors of global function fields. We also show improved lower bounds on
the corresponding function $\alpha_q^{\rm lin}$ for linear codes.
\endabstract

\vspace{1.0cm}

\noindent
\emph{Keywords:\/} Asymptotic theory of codes, Gilbert-Varshamov bound,
global function fields, Tsfasman-Vl\u{a}du\c{t}-Zink bound, Xing bound.

\vspace{0.5cm}

\noindent
2000 \emph{Mathematics Subject Classification:\/} Primary 11T71, 94B27,
94B65; Secondary 11R58, 14G50.

\vspace{0.8cm}

\newpage

\section{Introduction}
Let $\F_q$ be the finite field of order $q$, where $q$ is an arbitrary
prime power. For a  code $C$ over $\F_q$ (or in
other words a $q$-ary code), we denote by $n(C)$ its length and by $d(C)$
its minimum distance. We write $|M|$ for the cardinality of a finite set $M$.

For any prime power $q$, let $\alpha_q$ and $\alpha_q^{\rm lin}$ denote
the important functions in the asymptotic theory
of codes which are defined by
\be \label{definition.alpha.q.delta}
\alpha_q(\delta)={\rm sup} \, \{R \in [0,1]: (\delta,R) \in U_q \} \qquad
\mbox{for } 0 \le \delta \le 1
\ee
and
\be \label{definition.alpha.linear.q.delta}
\alpha_q^{\rm lin}(\delta) = {\rm sup} \, \{R \in [0,1]: (\delta, R) \in U_q^{\rm lin} \} \qquad \mbox{for } 0 \le \delta \le 1.
\ee
Here $U_q$ (resp. $U_q^{\rm lin}$) is the set of all ordered pairs $(\delta,R) \in [0,1]^2$ for which there exists a
sequence $\{C_i\}_{i=1}^\infty$ of not necessarily linear (resp. linear) codes over $\F_q$ such that $n(C_i) \ra \infty$
as $i \ra \infty$ and
\be
\delta=\lim_{i \ra \infty} \frac{d(C_i)}{n(C_i)}, \quad R=\lim_{i \ra \infty} \frac{\log_q |C_i|}{n(C_i)},
\nn\ee
where $\log_q$ is the logarithm to the base $q$. We refer to
\cite[Section 1.3.1]{TV} for some basic properties of the functions
$\alpha_q$ and $\alpha_q^{\rm lin}$. In particular, both functions are
nonincreasing on the interval $[0,1]$. Furthermore, we have the known
values $\alpha_q(0)=\alpha_q^{\rm lin}(0)=1$ and $\alpha_q(\delta)=
\alpha_q^{\rm lin}(\delta)=0$ for $(q-1)/q \le \delta \le 1$. It is trivial
that $\alpha_q(\delta) \ge \alpha_q^{\rm lin}(\delta)$ for $0 \le \delta
\le 1$.

A central problem in the asymptotic theory of codes is to find lower bounds
on $\alpha_q(\delta)$ for $0 < \delta < (q-1)/q$.
A classical lower bound is the asymptotic Gilbert-Varshamov bound which says that
{\small
\be \label{GV-bound}
\alpha_q^{\rm lin}(\delta) \ge R_{{\rm GV}}(\delta) := 1-\delta \log_q (q-1) + \delta
\log_q \delta + (1-\delta) \log_q (1-\delta)
\ee
}for $0 < \delta < (q-1)/q$. It is well known (see \cite[Section 6.2]{NX})
that for sufficiently large composite $q$ and for certain ranges
of the parameter $\delta$, one can beat the asymptotic
Gilbert-Varshamov bound by the Tsfasman-Vl\u{a}du\c{t}-Zink bound \cite{TVZ}
\be \label{TVZ-bound}
\alpha_q^{\rm lin}(\delta) \ge 1 - \delta - \frac{1}{A(q)} \qquad \mbox{for } 0 \le
\delta \le 1.
\ee
Here
$$
A(q) := \limsup_{g \ra \infty} \frac{N_q(g)}{g},
$$
where $N_q(g)$ denotes the maximum number of rational places that a global
function field of genus $g$ with full constant field $\F_q$ can have. We
recall from \cite[Chapter 5]{NX} that $A(q) > 0$ for all $q$ and that $A(q)=
\sqrt{q}-1$ if $q$ is a square. For nonsquares $q$ the exact value of $A(q)$
is not known, but we have lower and upper bounds on $A(q)$ (see again
\cite[Chapter 5]{NX}). We note, in
particular, the recent bound in \cite{BGS} which says that for any cube $q$ we
have
\be \label{eq14}
A(q) \ge \frac{2(q^{2/3}-1)}{q^{1/3}+2}.
\ee

The bound (\ref{TVZ-bound}) for $\alpha_q^{\rm lin}(\delta)$ was improved, although not uniformly in $\delta$,
by Vl\u{a}du\c{t} \cite{V} (see also \cite[Chapter 3.4]{TV}) and Xing
\cite{X1}. Elkies \cite{E} and Xing \cite{X2} considered not necessarily linear codes and Xing \cite{X2}
improved the bound (\ref{TVZ-bound}) for $\alpha_q(\delta)$ uniformly in $\delta$.
Shortly thereafter, Niederreiter and
\"Ozbudak \cite[Corollary 5.4]{NO1} improved the bound in Xing
\cite{X2} by showing that
\be \label{NO1-bound}
\alpha_q(\delta) \ge 1 - \delta - \frac{1}{A(q)} + \log_q \left(1+
\frac{1}{q^3} \right) \qquad \mbox{for } 0 \le \delta \le 1.
\ee
Later, Stichtenoth and Xing \cite{SX} gave a simpler proof of (\ref{NO1-bound}).

Recently, Niederreiter and \"Ozbudak \cite{NO2} improved the bound (\ref{NO1-bound}) for
certain values of $q$ and $\delta$.
In this paper we extensively refine and complement the methods of \cite{NO2}. We obtain further improvements on lower
bounds for $\alpha_q(\delta)$ and $\alpha^{\rm lin}_q(\delta)$ for certain values of $q$ and $\delta$ (see Theorem \ref{theorem.asymp}
and Corollary \ref{corollary.asym.linear}). In Section \ref{section.basic.construction} we present our basic code construction.
We obtain the cardinality of an important auxiliary set in this construction in Section \ref{section.cardinality}. Asymptotic
upper bounds on the cardinality of this set are given in Sections \ref{section.asymptotic.size.V.m.case.m.1} and \ref{section.asymptotic.size.V.m.case.m.general}.
We present our main results in Section \ref{section.asymptotic.constructon}. The final section is devoted to some
examples demonstrating the improvements obtained by the main results.

\section{The Basic Code Construction} \label{section.basic.construction}
In this section we present our basic construction of $q$-ary codes (see
Theorem \ref{theorem.basic} and Corollary \ref{corollary.basic.linear}).
We fix a global function field $F$ with full constant
field $\F_q$ and with at least one rational place.
Let $n \ge 1$ be the number of rational places of $F$ and let $P_1, \ldots, P_n$ be all rational places of $F$.
Let $h$ be the class number of $F$.
Let $v_P$ be the normalized discrete valuation of $F$ corresponding to the place $P$ of $F$.
Let $\P_F$ be the set of all places of $F$.
For $f \in F \setminus \{0\}$,
\be
(f)=\sum_{P \in \P_F} v_P(f) P
\nn\ee
denotes the principal divisor of $f$ and
\be
(f)_0=\sum_{\stackrel{P \in \P_F}{v_P(f) \ge 1}} v_P(f) P
\nn\ee
denotes the zero divisor of $f$.
For an arbitrary divisor
\be
G=\sum_{P \in \P_F} m_P P
\nn\ee
of $F$, we write $v_P(G)$ for the coefficient $m_P$ of $P$. We use the
standard notation
$$
\mathcal{L}(G)=\{f \in F: v_P(f) \ge -v_P(G) \ \mbox{for all } P \in \P_F\}
$$
for the Riemann-Roch space of $G$. In this section and in Section
\ref{section.cardinality}, all places and divisors are from the given
global function field $F$. We fix an integer $m \ge 1$.

\begin{definition} \label{definition.restrict.divisor}
For a positive divisor $D$, let $\overline{D}$ be the divisor
\be
\overline{D}=a_1P_1+ \cdots + a_nP_n,
\nn\ee
where $a_i=\min (m+1, v_{P_i}(D))$ for $1 \le i \le n$.
\end{definition}

\begin{definition} \label{definition.j.i.Jm}
For a positive divisor $D$, let
\be
\begin{array}{lcl}
j_0(D) & = & \left| \left\{i \in \{1, \ldots, n\}: v_{P_i}(D)=m \right\} \right|, \\
j_1(D) & = & \left| \left\{i \in \{1, \ldots, n\}: v_{P_i}(D)=m-1 \right\} \right|, \\
& \vdots & \\
j_m(D) & = & \left| \left\{i \in \{1, \ldots, n\}: v_{P_i}(D)=0 \right\} \right|. \\
\end{array}
\nn\ee
Moreover,  we  define
\be \label{relation.j.i.J.m}
J_m(D)=2j_1(D) + 3j_2(D) + \cdots + (m+1)j_m(D).
\ee
\end{definition}

\begin{definition} \label{definition.V.m}
For integers $r \ge s \ge 0$ and nonnegative integers $X_1,X_2, \ldots, X_m$, let
$\mathcal{V}_m(r,s;X_1,X_2, \ldots, X_m)$ be the
set consisting of the positive divisors $D$ of the global function field $F$
satisfying all of the following:
\begin{itemize}

\item Condition 1: $\deg (D)=r$ and $\deg \left(\overline{D}\right) \ge s$,

\item Condition 2:
\be
\begin{array}{rcl}
j_m(D) &\le & 2X_m, \\
j_{m-1}(D) &\le & 2 X_{m-1} + X_m,\\
j_{m-2}(D) & \le & 2 X_{m-2} + \left(X_{m-1}+X_m\right), \\
&\vdots& \\
j_1(D) & \le & 2X_1+ \left(X_2+X_3+ \cdots + X_m\right),
\end{array}
\nn\ee

\item Condition 3: $J_m(D) \le 2\left(2X_1+3X_2+ \cdots + (m+1)X_m \right)$.

\end{itemize}
\end{definition}

\begin{proposition} \label{proposition.existence.G}
For integers $r\ge s \ge 0$ and nonnegative integers $X_1, \ldots, X_m$, if
\be
\left|\mathcal{V}_m(r,s;X_1, \ldots, X_m)\right| < h,
\nn\ee
then there exists a
divisor $G$ of degree $r$ such that $\supp (G) \cap \{P_1, \ldots, P_n\}=\emptyset$ and
for each $f \in \mathcal{L}(G) \setminus \{0\}$,
if $E=(f)_0$ satisfies  Conditions 2 and 3 of Definition \ref{definition.V.m}
with the given $X_1, \ldots, X_m$,
then $\deg \left( \overline{E} \right) \le s-1$.
\end{proposition}
\begin{proof}
As $\left| \mathcal{V}_m(r,s;X_1, \ldots, X_m) \right| < h$, there exists a degree
$r$ divisor $G$ such that $G \not \sim V$ for any $V \in
\mathcal{V}_m(r,s;X_1, \ldots, X_m)$. Using the Weak Approximation Theorem
\cite[Theorem I.3.1]{S}, we can assume that $\supp (G) \cap \{P_1,
\ldots, P_n\}=\emptyset$ without loss of generality (compare with \cite[Proof
of Corollary 2.2]{NO2}). Let $f \in \mathcal{L}(G) \setminus \{0\}$,
$D=G+ (f)$, and $E=(f)_0$. Since
$\supp (G) \cap \{P_1, \ldots, P_n\}=\emptyset$ and $D$ is
positive, we have $\overline{D}=\overline{E}$. Assume that
Conditions 2 and 3 of Definition \ref{definition.V.m} are satisfied by
$E$. If  $\deg \left(\overline{E}\right) \ge s$, then $D \in
\mathcal{V}_m(r,s;X_1, \ldots, X_m)$ and hence $D \not \sim G$, which is a contradiction. Thus, we must have $\deg \left(\overline{E}\right) \le s-1$.
\end{proof}

Now give another definition related to our construction.

\begin{definition} \label{definition.I.i}
For ${\boldsymbol \alpha}=\left(\alpha_1^{(1)}, \ldots, \alpha_m^{(1)}, \alpha_1^{(2)}, \ldots, \alpha_m^{(2)},\ldots \ldots, \alpha_1^{(n)}, \ldots, \alpha_m^{(n)} \right)
\in \F_q^{mn}$, let $I_m({\boldsymbol \alpha})$, $I_{m-1}({\boldsymbol \alpha}), \ldots, I_1({\boldsymbol \alpha})$ be the subsets of
$\{1, \ldots, n\}$ defined by
\be
\begin{array}{lcl}
\dd I_m({\boldsymbol \alpha}) & = & \dd \left\{ i \in \{1, \ldots, n\}: \alpha_m^{(i)} \neq 0 \right\}, \\
\dd I_{m-1}({\boldsymbol \alpha}) & = & \dd  \left\{ i \in \{1, \ldots, n\}: \alpha_m^{(i)}=0, \; \alpha_{m-1}^{(i)} \neq 0 \right\}, \\
& \vdots & \\
\dd I_1({\boldsymbol \alpha}) & = & \dd \left\{i \in \{1, \ldots, n\}: \alpha_m^{(i)}=\cdots=\alpha_2^{(i)}=0, \; \alpha_1^{(i)} \neq 0\right\}.
\end{array}
\nn\ee
\end{definition}

The following two lemmas are related to Definition \ref{definition.I.i} and important for our construction.

\begin{lemma} \label{lemma.bound.on.I.wighted.cardinality.sum}
For $\boldsymbol{\alpha}, \boldsymbol{\beta} \in \F_q^{mn}$, we have
\be
\begin{array}{l}
\dd 2 \left|I_1(\boldsymbol{\alpha} - \boldsymbol{\beta}) \right| + 3 \left| I_2(\boldsymbol{\alpha}-\boldsymbol{\beta}) \right| + \cdots +
(m+1) \left| I_m(\boldsymbol{\alpha} - \boldsymbol{\beta}) \right| \\ \\
\dd \le
2 \left|I_1(\boldsymbol{\alpha} ) \right| + 3 \left| I_2(\boldsymbol{\alpha}) \right| + \cdots +
(m+1) \left| I_m(\boldsymbol{\alpha} ) \right| \\ \\
\dd
\hspace{0.5cm} +2 \left|I_1(\boldsymbol{\beta} ) \right| + 3 \left| I_2(\boldsymbol{\beta}) \right| + \cdots +
(m+1) \left| I_m(\boldsymbol{\beta} ) \right|.
\end{array}
\nn\ee
\end{lemma}
\begin{proof}
Let
$\boldsymbol{\alpha}=\big(\alpha_1^{(1)}, \ldots, \alpha_m^{(1)},\ldots \ldots, \alpha_1^{(n)}, \ldots, \alpha_m^{(n)} \big)$
and
$\boldsymbol{\beta}=\big(\beta_1^{(1)}, \ldots, \beta_m^{(1)},\ldots$ $ \ldots, \beta_1^{(n)}, \ldots, \beta_m^{(n)} \big)$.
Let $A \subseteq \{1, \ldots, n\}$ be the set consisting of the $i \in \{1, \ldots, n\}$ such that
$\left(\alpha_1^{(i)}, \ldots, \alpha_m^{(i)}\right) \neq \boldsymbol{0}$ or
$\left(\beta_1^{(i)}, \ldots, \beta_m^{(i)}\right) \neq \boldsymbol{0}$.
If $A=\emptyset$, then $\boldsymbol{\alpha}=\boldsymbol{\beta} = \boldsymbol{\alpha} - \boldsymbol{\beta}=\boldsymbol{0}$
and the result follows immediately. If $A \neq \emptyset$, then for each $i \in A$, let $1 \le \ell_i \le m$ be the largest integer such that
$\alpha_{\ell_i}^{(i)} \neq 0$ or $\beta_{\ell_i}^{(i)} \neq 0$. For each $i \in A$, we have
\be
i \not \in \bigcup_{\ell_i < j \le m} I_j(\boldsymbol{\alpha} - \boldsymbol{\beta}),
\nn\ee
and also $i \in I_{\ell_i}(\boldsymbol{\alpha})$ or $i \in I_{\ell_i}(\boldsymbol{\beta})$. Hence for each $i \in A$
we obtain
\be
\begin{array}{l}
\dd 2 \left|\{i\} \cap I_1(\boldsymbol{\alpha} - \boldsymbol{\beta}) \right| + 3 \left| \{i\} \cap I_2(\boldsymbol{\alpha}-\boldsymbol{\beta}) \right| + \cdots +
(m+1) \left|\{i\} \cap I_m(\boldsymbol{\alpha} - \boldsymbol{\beta}) \right| \\ \\
\dd \le
2 \left|\{i\} \cap I_1(\boldsymbol{\alpha} ) \right| + 3 \left| \{i\} \cap I_2(\boldsymbol{\alpha}) \right| + \cdots +
(m+1) \left| \{i\} \cap I_m(\boldsymbol{\alpha} ) \right| \\ \\
\dd
\hspace{0.5cm} +2 \left|\{i\} \cap I_1(\boldsymbol{\beta} ) \right| + 3 \left|\{i\} \cap I_2(\boldsymbol{\beta}) \right| + \cdots +
(m+1) \left| \{i\} \cap I_m(\boldsymbol{\beta} ) \right|.
\end{array}
\nn\ee
We complete the proof by summing over all $i \in A$.
\end{proof}

\begin{lemma} \label{lemma.second.bound.on.I.i}
For $\boldsymbol{\alpha}, \boldsymbol{\beta} \in \F_q^{mn}$, we have the following
containment relations:
\be
\begin{array}{lcl}
I_{m}(\boldsymbol{\alpha} - \boldsymbol{\beta}) & \subseteq & I_{m}(\boldsymbol{\alpha}) \cup I_m(\boldsymbol{\beta}), \\
I_{m-1}(\boldsymbol{\alpha} - \boldsymbol{\beta}) & \subseteq & I_{m-1}(\boldsymbol{\alpha}) \cup
I_{m-1}(\boldsymbol{\beta}) \cup \left\{ I_m(\boldsymbol{\alpha}) \cap I_m(\boldsymbol{\beta}) \right\}, \\
I_{m-2}(\boldsymbol{\alpha} - \boldsymbol{\beta}) & \subseteq & I_{m-2}(\boldsymbol{\alpha}) \cup I_{m-2}(\boldsymbol{\beta})
\cup \left\{ I_{m-1}(\boldsymbol{\alpha}) \cap I_{m-1}(\boldsymbol{\beta}) \right\} \cup \{I_m(\boldsymbol{\alpha}) \cap I_m(\boldsymbol{\beta}) \},\\
& \vdots & \\
I_1(\boldsymbol{\alpha} - \boldsymbol{\beta}) & \subseteq &
I_1(\boldsymbol{\alpha}) \cup I_1(\boldsymbol{\beta}) \cup \bigcup_{2 \le \nu \le m} \left\{ I_\nu(\boldsymbol{\alpha}) \cap I_\nu(\boldsymbol{\beta}) \right\}.
\end{array}
\nn\ee
\end{lemma}
\begin{proof}
First we consider the case of the subscript $m$ and we assume that $i \in I_m(\boldsymbol{\alpha} -\boldsymbol{\beta})$.
Then $\alpha_m^{(i)} \neq \beta_m^{(i)}$ and at least one of $\alpha_m^{(i)}$ and $\beta_m^{(i)}$ is nonzero.
Hence $i \in I_m(\boldsymbol{\alpha}) \cup I_m(\boldsymbol{\beta})$.

Next we consider the case of the subscript $m-1$ and we assume that $i \in I_{m-1}(\boldsymbol{\alpha}-\boldsymbol{\beta})$.
We have $\alpha_m^{(i)}=\beta_m^{(i)}$ and $\alpha_{m-1}^{(i)} \neq \beta_{m-1}^{(i)}$. If $\alpha_m^{(i)}=\beta_m^{(i)} \neq 0$,
then $i \in I_m(\boldsymbol{\alpha}) \cap I_m(\boldsymbol{\beta})$. If $\alpha_m^{(i)}=\beta_m^{(i)} = 0$, then since
at least one of $\alpha_{m-1}^{(i)}$ and $\beta_{m-1}^{(i)}$ is nonzero, we get $i \in I_{m-1}(\boldsymbol{\alpha}) \cup I_{m-1}(\boldsymbol{\beta})$.

Now we consider the case of  the subscript $m-2$.
Assume that
$i \in I_{m-2}(\boldsymbol{\alpha}- \boldsymbol{\beta})$. Then $\alpha_m^{(i)}=\beta_m^{(i)}$,
$\alpha_{m-1}^{(i)}=\beta_{m-1}^{(i)}$, and $\alpha_{m-2}^{(i)} \neq \beta_{m-2}^{(i)}$. If
$\alpha_m^{(i)}=\beta_m^{(i)} \neq 0$,
then $i \in I_m(\boldsymbol{\alpha}) \cap I_m(\boldsymbol{\beta})$. If $\alpha_m^{(i)}=\beta_m^{(i)}=0$
and $\alpha_{m-1}^{(i)}=\beta_{m-1}^{(i)} \neq 0$, then $i \in I_{m-1}(\boldsymbol{\alpha}) \cap I_{m-1}(\boldsymbol{\beta})$.
Finally, if $\alpha_m^{(i)}=\beta_m^{(i)}=0$
and $\alpha_{m-1}^{(i)}=\beta_{m-1}^{(i)} = 0$, then since $\alpha_{m-2}^{(i)}$ and $\beta_{m-2}^{(i)}$ are distinct, we get
$i \in I_{m-2}(\boldsymbol{\alpha})$ or $i \in I_{m-2}(\boldsymbol{\beta})$. We complete the proof similarly for each subscript
$1 \le \nu \le m$.
\end{proof}

For each $i=1, \ldots, n$, let $t_i$ be a local parameter of $F$ at $P_i$.
Assume that $G$ is a divisor with $\supp(G) \cap \{P_1, \ldots, P_n\}=\emptyset$
and $\dim \left(\mathcal{L}(G)\right) \ge 1$.
For $f$ in the Riemann-Roch space $\mathcal{L}(G)$, the local expansion of $f$ at $P_i$ has the form
\be
f=\sum_{l=0}^\infty f^{(l)}(P_i)t_i^l
\nn\ee
with $f^{(l)}(P_i) \in \F_q$ for $1 \le i \le n$ and $l \ge 0$.
For each $i=1, \ldots, n$, let
\be
\begin{array}{rcl}
\phi_i: \mathcal{L}(G) & \ra & \F_q^m \\
f & \mapsto & \left( f^{(m-1)}(P_i), \ldots, f^{(1)}(P_i), f^{(0)}(P_i) \right).
\end{array}
\nn\ee
Let $\boldsymbol{\Phi}$ be the $\F_q$-linear map defined by
\be \label{definition.Phi}
\begin{array}{rcl}
\Phi:\mathcal{L}(G) & \ra & \F_q^{mn} \\
f & \mapsto & \left( \phi_1(f), \ldots, \phi_n(f) \right).
\end{array}
\ee
Moreover, let $\psi$ be the $\F_q$-linear map
\be \label{definition.psi}
\begin{array}{rcl}
\psi:\mathcal{L}(G) & \ra & \F_q^{n} \\
f & \mapsto & \left( f^{(m)}(P_1), \ldots, f^{(m)}(P_n) \right).
\end{array}
\ee

\begin{lemma} \label{lemma.relation.J.m.I.i}
For a divisor $G$ with $\supp(G) \cap \{P_1, \ldots, P_n\}=\emptyset$ and $\dim \left(\mathcal{L}(G)\right) \ge 1$, let
$f \in \mathcal{L}(G) \setminus \{0\}$. Moreover, let $E=(f)_0$ be the zero divisor of $f$ and
$\boldsymbol{\alpha}:=\boldsymbol{\Phi}(f) \in \F_q^{mn}$. Then
\be
j_1(E)=|I_1(\boldsymbol{\alpha})|, \; j_2(E)=|I_2(\boldsymbol{\alpha})|,  \; \ldots \; , j_m(E)=|I_m(\boldsymbol{\alpha})|,
\nn\ee
and
\be
J_m(E)=2 \left|I_1(\boldsymbol{\alpha})\right| + 3 \left|I_2(\boldsymbol{\alpha})\right| + \cdots + (m+1)  \left|I_m(\boldsymbol{\alpha})\right|.
\nn\ee
\end{lemma}
\begin{proof}
For each $1 \le i \le n$,
using Definition \ref{definition.I.i} we obtain
\be
\begin{array}{lcl}
i \in I_m(\boldsymbol{\alpha}) & \iff & v_{P_i}(E)=0, \\
i \in I_{m-1}(\boldsymbol{\alpha}) & \iff & v_{P_i}(E)=1, \\
 & \vdots &  \\
 i \in I_1(\boldsymbol{\alpha}) &\iff &  v_{P_i}(E)=m-1.
\end{array}
\nn\ee
Hence by Definition \ref{definition.j.i.Jm} we have
\be
j_m(E)=|I_m(\boldsymbol{\alpha})|, \; j_{m-1}(E)=|I_{m-1}(\boldsymbol{\alpha})|,  \; \ldots \; , j_1(E)=|I_1(\boldsymbol{\alpha})|.
\nn\ee
Using (\ref{relation.j.i.J.m}) we complete the proof.
\end{proof}

For ${\boldsymbol c} \in \F_q^{mn}$ and nonnegative real numbers $x_1, \ldots, x_m$ with
$x_1+ \cdots + x_m  \le 1$, let $M(x_1, \ldots, x_m;{\boldsymbol c})$ be the subset of $\F_q^{mn}$ defined by
\be
M(x_1, \ldots, x_m; \boldsymbol{c})=\left\{\boldsymbol{\alpha} \in \F_q^{mn}:
\left| I_1(\boldsymbol{\alpha}-\boldsymbol{c})\right| \le \lfloor x_1 n \rfloor, \ldots,
\left| I_m(\boldsymbol{\alpha}-\boldsymbol{c}) \right| \le \lfloor x_m n \rfloor \right\}.
\nn\ee
We have
\be \label{eq1.basic.const}
\begin{array}{l}
\dd \left| M(x_1, \ldots, x_m; \boldsymbol{c}) \right| = \left| M(x_1, \ldots, x_m; \boldsymbol{0}) \right| \\ \\
\ge \left| \left\{ \boldsymbol{\alpha} \in \F_{q^{mn}}: \left| I_1(\boldsymbol{\alpha}) \right|=\lfloor x_1 n \rfloor, \ldots,
\left|I_m(\boldsymbol{\alpha}) \right| = \lfloor x_m n \rfloor \right\} \right| \\ \\
\dd = {n \choose \lfloor x_m n \rfloor } (q-1)^{\lfloor x_m n \rfloor} q^{(m-1)\lfloor x_m n \rfloor} \\ \\
\dd \times {n-\lfloor x_m n \rfloor \choose \lfloor x_{m-1} n \rfloor } (q-1)^{\lfloor x_{m-1} n \rfloor} q^{(m-2)\lfloor x_{m-1} n \rfloor} \\ \\
\dd \times \cdots \\ \\
\dd \times { n-\left( \lfloor x_m n \rfloor + \lfloor x_{m-1} n \rfloor + \cdots + \lfloor x_2 n \rfloor \right) \choose \lfloor x_1 n \rfloor } (q-1)^{\lfloor x_1 n \rfloor}.
\end{array}
\ee

Now we are ready to give our basic code construction. Assume that $r \ge s \ge 0$ are integers and $ x_1, \ldots, x_m \ge 0$
are real numbers such that
\be \label{assumption1}
\left|\mathcal{V}_m \left(r,s;\lfloor x_1 n \rfloor,  \lfloor x_2 n \rfloor,  \ldots ,  \lfloor x_m n \rfloor  \right) \right|< h.
\ee
Let $G$ be a divisor of degree $r$ obtained using (\ref{assumption1}) and Proposition \ref{proposition.existence.G}.
Recall the linear maps $\boldsymbol{\Phi}$ and $\psi$ defined in (\ref{definition.Phi}) and (\ref{definition.psi}),
respectively, using the chosen divisor $G$.
The map $\boldsymbol{\Phi}$ is not necessarily surjective. If
\be \label{assumption2}
\left| \mathcal{L}(G) \right| \cdot \left| M(x_1, \ldots, x_m; \boldsymbol{0}) \right| > q^{mn},
\ee
then there exists $\boldsymbol{c} \in \F_q^{mn}$ such that for the set
\be \label{definition.N.c}
N_{\boldsymbol{c}}:=\left\{ f \in \mathcal{L}(G) : \boldsymbol{\Phi}(f) \in M(x_1, \ldots, x_m; \boldsymbol{c}) \right\}
\ee
we have
\be \label{N.c.lower.bound}
\left| N_{\boldsymbol{c}} \right| \ge \frac{\left| \mathcal{L}(G) \right| \cdot \left| M(x_1, \ldots, x_m; \boldsymbol{0}) \right|}{q^{mn}} > 1.
\ee

\begin{theorem} \label{theorem.basic}
Assume that $r \ge s \ge 0$ are integers and that $x_1, \ldots, x_m$
are nonnegative real numbers with $x_1+ \cdots + x_m \le 1$ satisfying (\ref{assumption1}).
Let $G$ be a divisor of degree $r$ obtained using (\ref{assumption1}) and Proposition \ref{proposition.existence.G}.
Assume also that (\ref{assumption2}) holds and that
\be \label{e1.theorem.basic}
(m+1)n \ge s +2 \sum_{l=1}^m (l+1) \lfloor x_l n \rfloor.
\ee
Using the chosen divisor $G$ and (\ref{assumption2}),
let $\boldsymbol{c} \in \F_q^{mn}$ be
such that the set $N_{\boldsymbol{c}}$ satisfies (\ref{N.c.lower.bound}). Let $C$
be the $q$-ary code of length $n$ given by $C =\psi\left(N_{\boldsymbol{c}}\right)$.
Then for the cardinality $|C|$ of $C$ we have
\be
|C| \ge \left\lceil \frac{ \mathcal{L}(G) \cdot \left| M(x_1, \ldots, x_m; \boldsymbol{0}) \right|}{q^{mn}} \right\rceil
\nn\ee
and for the minimum distance $d(C)$ of $C$ we have
\be
d(C) \ge (m+1)n +1 -s - 2 \sum_{l=1}^m (l+1) \lfloor x_l n \rfloor.
\nn\ee
\end{theorem}
\begin{proof}
Let $f_1,f_2 \in N_{\boldsymbol{c}}$ be such that $f_1 \neq f_2$ and put $f=f_1-f_2 \in \mathcal{L}(G)$. Let $E$ be the
zero divisor of $f$ and
\be
\overline{E}=a_1P_1+ \cdots + a_n P_n
\nn\ee
be the divisor defined in Definition \ref{definition.restrict.divisor}. Let $\boldsymbol{\Phi}(f_1)=\boldsymbol{\alpha}$
and $\boldsymbol{\Phi}(f_2)=\boldsymbol{\beta}$. We have
\be \label{ep1.theorem.basic}
\boldsymbol{\Phi}(f) =\boldsymbol{\alpha}-\boldsymbol{\beta}.
\ee
As $\boldsymbol{\alpha}, \boldsymbol{\beta} \in M(x_1, \ldots, x_m;\boldsymbol{c})$, we also have
\be \label{ep2.theorem.basic}
\left| I_i(\boldsymbol{\alpha}-\boldsymbol{c})\right| \le \lfloor x_i n\rfloor \ \mbox{and} \
\left| I_i(\boldsymbol{\beta} - \boldsymbol{c})\right| \le \lfloor x_i n \rfloor \; \mbox{for $1 \le i \le n$}.
\ee
Using (\ref{ep1.theorem.basic}), (\ref{ep2.theorem.basic}), Lemmas \ref{lemma.relation.J.m.I.i} and \ref{lemma.bound.on.I.wighted.cardinality.sum},
we obtain that
\be
\begin{array}{rl}
\dd J_m(E)= & 2 \left|I_1(\boldsymbol{\alpha}-\boldsymbol{\beta})\right|
+3 \left|I_2(\boldsymbol{\alpha}-\boldsymbol{\beta})\right|
+\cdots
+(m+1) \left|I_m(\boldsymbol{\alpha}-\boldsymbol{\beta})\right| \\ \\

\le & \dd 2 \left|I_1(\boldsymbol{\alpha}-\boldsymbol{c})\right|
+3 \left|I_2(\boldsymbol{\alpha}-\boldsymbol{c})\right|
+\cdots
+(m+1) \left|I_m(\boldsymbol{\alpha}-\boldsymbol{c})\right| \\ \\
& \dd \hspace{0.1cm} +  2 \left|I_1(\boldsymbol{\beta}-\boldsymbol{c})\right|
+3 \left|I_2(\boldsymbol{\beta}-\boldsymbol{c})\right|
+\cdots
+(m+1) \left|I_m(\boldsymbol{\beta}-\boldsymbol{c})\right| \\ \\
\dd \le & \dd 2 \left( 2 \lfloor x_1 n \rfloor + 3 \lfloor x_2 n \rfloor + \cdots + (m+1) \lfloor x_m n \rfloor \right).
\end{array}
\nn\ee
Moreover, using (\ref{ep1.theorem.basic}), (\ref{ep2.theorem.basic}), Lemmas \ref{lemma.relation.J.m.I.i} and \ref{lemma.second.bound.on.I.i},
we further obtain that
\be
\begin{array}{lcl}
j_m(E) & = & |I_m((\boldsymbol{\alpha}-\boldsymbol{c})-(\boldsymbol{\beta}-\boldsymbol{c}))|
 \le   |I_m(\boldsymbol{\alpha}-\boldsymbol{c})| + |I_m(\boldsymbol{\beta}-\boldsymbol{c})| \le 2 \lfloor x_m n \rfloor, \\ \\
j_{m-1}(E) & = & |I_{m-1}((\boldsymbol{\alpha}-\boldsymbol{c})-(\boldsymbol{\beta}-\boldsymbol{c}))| \\
& \le&  |I_{m-1}(\boldsymbol{\alpha}-\boldsymbol{c})|  + |I_{m-1}(\boldsymbol{\beta}-\boldsymbol{c})|
+ |I_m(\boldsymbol{\alpha}-\boldsymbol{c}) \cap I_m(\boldsymbol{\beta}-\boldsymbol{c})| \\
&  \le & 2  \lfloor x_{m-1}n \rfloor + \lfloor x_m n \rfloor, \\ \\
& \vdots & \\
j_1(E) & = & |I_1((\boldsymbol{\alpha}-\boldsymbol{c})-(\boldsymbol{\beta}-\boldsymbol{c}))| \\
& \le & |I_1(\boldsymbol{\alpha}-\boldsymbol{c})| + |I_1(\boldsymbol{\beta}-\boldsymbol{c})| + \sum_{\nu=2}^m |I_\nu(\boldsymbol{\alpha}-\boldsymbol{c}) \cap I_\nu(\boldsymbol{\beta}-\boldsymbol{c})| \\
& \le & 2\lfloor x_1 n \rfloor + \sum_{\nu=2}^m \lfloor x_\nu n \rfloor.
\end{array}
\nn\ee
Hence by the choice of the divisor $G$ (cf. Proposition \ref{proposition.existence.G}), we have
\be \label{ep3.theorem.basic}
\deg \left(\overline{E}\right) \le s-1.
\ee
Moreover, we obtain
\be
\sum_{i=1}^n (m+1-a_i)=(m+1)n - \sum_{i=1}^n a_i = (m+1)n - \deg \left(\overline{E}\right) \ge (m+1)n -s +1,
\nn\ee
where we used (\ref{ep3.theorem.basic}). Let $||\psi(f)||$ denote the Hamming weight of the vector $\psi(f) \in \F_q^n$.
Then using Definition \ref{definition.j.i.Jm} and (\ref{relation.j.i.J.m}), we have
\be
\begin{array}{rl}
\dd \sum_{i=1}^n (m+1-a_i)= & \dd \sum_{\stackrel{i=1}{0 \le a_i \le m}}^n (m+1-a_i) \le
||\psi(f)|| + \sum_{\stackrel{i=1}{0 \le a_i \le m-1}}^n (m+1-a_i) \\ \\
= & \dd ||\psi(f)|| + J_m(E).
\end{array}
\nn\ee
Therefore we obtain
\be
\begin{array}{rl}
\dd ||\psi(f)|| & \ge (m+1)n-s+1-J_m(E) \\
&\dd \ge (m+1)n-s+1-
2 \left( 2 \lfloor x_1 n \rfloor + 3 \lfloor x_2 n \rfloor + \cdots + (m+1) \lfloor x_m n \rfloor \right).
\end{array}
\nn\ee
Using (\ref{e1.theorem.basic}) we obtain that $d(C) \ge 1$, and so the map $\psi$ is one-to-one on $N_{\boldsymbol{c}}$.
Therefore $|C|= |N_{\boldsymbol{c}}|$, and hence the lower bound on $|C|$ follows from (\ref{N.c.lower.bound}). This completes the proof.
\end{proof}

In a special case related to Theorem \ref{theorem.basic}, we make sure to construct linear codes. Later in this paper,
the following result will be used to obtain lower bounds on the function $\alpha_q^{\rm lin}(\delta)$,  which is defined in (\ref{definition.alpha.linear.q.delta}).

\begin{corollary} \label{corollary.basic.linear}
Assume that $r \ge s \ge 0$ are integers and that $x_1=x_2=\cdots=x_m=0$  satisfy (\ref{assumption1}).
Let $G$ be a divisor of degree $r$ obtained using (\ref{assumption1}) and Proposition \ref{proposition.existence.G}.
Assume also that
\be \label{e1.corollary.basic.linear}
\left|\mathcal{L}(G)\right| > q^{mn}
\ee
and that $(m+1)n \ge s$. Using the chosen divisor $G$ and the kernel  of the corresponding map $\boldsymbol{\Phi}$, put
$C=\psi\left( \Ker \; \boldsymbol{\Phi} \right)$.
Then $C$ is a linear code over $\F_q$ of length $n$. Moreover,
for the dimension of $C$ we have
\be
\dim (C) \ge \dim \left( \mathcal{L}(G)\right) -mn
\nn\ee
and for the minimum distance $d(C)$ of $C$ we have
\be
d(C) \ge (m+1)n +1 -s.
\nn\ee
\end{corollary}
\begin{proof}
The kernel  of $\boldsymbol{\Phi}$ is an $\F_q$-linear subspace of $\mathcal{L}(G)$ and is the
Riemann-Roch space given by
\be
\Ker \; \boldsymbol{\Phi} =\mathcal{L}\left(G-m(P_1+ \cdots+P_n) \right).
\nn\ee
As $\dim \left(\mathcal{L}\left(G-m(P_1+ \cdots+P_n) \right)\right) \ge \dim \left(\mathcal{L}(G)\right) - mn$, using (\ref{e1.corollary.basic.linear})
we obtain that $\Ker \; \boldsymbol{\Phi}  \neq \{0\}$. The maps $\boldsymbol{\Phi}$ and $\psi$ are
$\F_q$-linear, and hence $C$ is a linear code over $\F_q$. We obtain the bounds on the dimension and the minimum distance of $C$
using similar methods as in the proof of Theorem \ref{theorem.basic}.
\end{proof}
\begin{remark} \label{remark.on.corollary.basic.linear}
For $x_1=x_2=\cdots=x_m=0$, the conditions (\ref{assumption2}) and (\ref{e1.corollary.basic.linear}) are equivalent.
\end{remark}

\section{The Cardinality of $\mathcal{V}_m(r,s;X_1, \ldots,X_m)$} \label{section.cardinality}

In this section we will compute the cardinality of the set
$\mathcal{V}_m(r,s;X_1, \ldots,X_m)$ for integers $r \ge s \ge 0$ and nonnegative integers $X_1, \ldots, X_m$
(see Definition \ref{definition.V.m} for the definition of this set).
The notation we introduced in Section \ref{section.basic.construction}
remains operative.

\begin{lemma} \label{lemma.relation.on.j.and.degree}
For any positive divisor $D$, we have
\be
\deg \left(\overline{D}\right) + j_0(D) + 2j_1(D) + \cdots + (m+1)j_m(D)=(m+1)n.
\nn\ee
\end{lemma}
\begin{proof}
For $0 \le \ell \le m$, let $S_\ell=\left\{P \in \{P_1, \ldots, P_n\}: v_{P}(\overline{D})=m-\ell\right\}$. Note that
$|S_\ell|=j_\ell(D)$ for each $0 \le \ell \le m$. We have
\be
\sum_{P \in \{P_1, \ldots, P_n\}} (m+1-v_P(\overline{D}))=(m+1)n-\deg \left(\overline{D}\right)
\nn\ee
and also
\be
\begin{array}{rl}
\dd \sum_{P \in \{P_1, \ldots, P_n\}} \left( m+1-v_{P}(\overline{D}) \right) & = \dd \sum_{\ell=0}^m \sum_{P \in S_\ell} \left( m+1 - v_P(\overline{D}) \right) \\ \\
& =\dd \sum_{\ell=0}^m \sum_{P \in S_\ell} (\ell+1)
 = \dd \sum_{\ell=0}^m (\ell+1) j_\ell(D).
\end{array}
\nn\ee
This completes the proof.
\end{proof}

\begin{definition} \label{definition.set.U}
For integers $r \ge t\ge 0$ and $j_1, \ldots, j_m \ge 0$, let $\mathcal{U}(r,t;j_1, \ldots, j_m)$ be the set
of positive divisors given by
\be
\mathcal{U}(r,t;j_1, \ldots, j_m)=\left\{D\ge 0: \deg \left(D\right)=r, \; \deg \left(\overline{D}\right)=t, \; j_1(D)=j_1, \ldots, j_m(D)=j_m\right\}.
\nn\ee
\end{definition}

\begin{lemma} \label{lemma.non.emptiness.mathcal.u}
For integers $r \ge t \ge 0$ and $j_1, \ldots, j_m\ge 0$, the set $\mathcal{U}(r,t;j_1, \ldots, j_m)$ is not empty
if and only if
\be
mn-\left(j_1+2j_2+ \cdots +m j_m\right) \le t \le (m+1)n-\left(2j_1+3j_2+ \cdots + (m+1)j_m\right)
\nn\ee
holds and also provided that
there exists a degree $r-t$ positive divisor whose support is disjoint from the set $\{P_1, \ldots, P_n\}$
when $mn=t+j_1+ 2j_2+ \cdots +mj_m$ and $r>t$.
\end{lemma}
\begin{proof}
Let $D \in \mathcal{U}(r,t;j_1, \ldots, j_m)$. Using Lemma \ref{lemma.relation.on.j.and.degree} we have
\be \label{ep0.lemma.non.emptiness.mathcal.u}
j_0(D)=(m+1)n-\left( 2j_1(D)+ \cdots + (m+1)j_m(D) \right)-t,
\ee
and so in particular
\be
t \le (m+1)n - \left(2 j_1+ 3j_2+ \cdots + (m+1) j_m\right).
\nn\ee
Moreover by definition of $\overline{D}$,
\be
\begin{array}{rl}
t \ge & j_{m-1}(D) + 2j_{m-2}(D) + \cdots + mj_0(D) \\
= & j_{m-1}(D) + 2j_{m-2}(D) + \cdots + (m-1)j_1(D) \\
& + m(m+1)n - \left( 2mj_1(D) + \cdots + (m+1)mj_m(D) \right)  -mt,
\end{array}
\nn\ee
where we used (\ref{ep0.lemma.non.emptiness.mathcal.u}) in the second step.
Therefore
\be
\begin{array}{rl}
\dd (m+1)t \ge & \dd (m+1)mn  \\ \\
& \dd - \big( (m+1)mj_m(D) + (m^2-1) j_{m-1}(D) + \left((m-1)m-2\right)j_{m-2}(D) \\ \\
& \dd \hspace{0.6cm}+ \cdots + \left(2m-(m-1)\right)j_1(D) \big) \\ \\
= & \dd (m+1)mn - (m+1) \left( mj_m(D) + (m-1)j_{m-1}(D) + \cdots + j_1(D)\right),
\end{array}
\nn\ee
which means that
\be \label{ep1.lemma.non.emptiness.mathcal.u}
t \ge mn - \left(j_1 + 2j_2 + \cdots + mj_m \right).
\ee
Also, if this is an equality, then the set $\left\{P \in \{P_1, \ldots, P_n\}: v_P(D) \ge m+1\right\}$ is empty.
Therefore, if  equality in (\ref{ep1.lemma.non.emptiness.mathcal.u}) holds and $r>t$, then there exists
a positive divisor
of degree $r-t$ whose support is disjoint from $\{P_1, \ldots, P_n\}$.

Now we prove the converse. Let $S_m=\{1, \ldots, j_m\}$, $S_{m-1}=\{j_m+1, \ldots, j_m+j_{m-1}\}$,
\ldots, $S_1=\{(j_m+ \cdots+j_2)+1, \ldots, (j_m+ \cdots + j_2)+j_1\}$.
They are pairwise disjoint sets of natural numbers.
We note that for each $1 \le \ell \le m$, we have
$|S_\ell|=j_{\ell}$. Comparing both sides of the inequalities for
$t$ given in the statement of the lemma, we obtain that
\be
j_1+j_2+ \cdots + j_m \le n.
\nn\ee
Let
\be \label{ep2.lemma.non.emptiness.mathcal.u}
j_0=(m+1)n-\left(2j_1+3j_2+ \cdots + (m+1)j_m\right)-t.
\ee
Using the upper bound on
$t$ in the statement of the lemma, we get $j_0 \ge 0$. Moreover, using
$t \ge mn-(j_1+2j_2+ \cdots+ mj_m)$
we obtain
\be
j_0+j_1+\cdots + j_m=(m+1)n-\left(j_1+2j_2+ \cdots + mj_m\right)-t \le n.
\nn\ee
Let $S_0=\{(j_m+ \cdots + j_1)+1, \ldots, (j_m+ \cdots + j_1)+j_0\}$. Note that $S_0, \ldots, S_m$ are pairwise disjoint
subsets of $\{1, \ldots, n\}$. For each $i \in \{1, \ldots,n\}$, let
\be
a_i=\left\{
\begin{array}{ll}
m-\ell & \mbox{if $i \in S_\ell$ for some $0 \le \ell \le m$}, \\
m+1 & \mbox{otherwise}.
\end{array}
\right.
\nn\ee
Assume that $j_m+ \cdots + j_1+j_0 < n$ and put
\be
D=(r-t)P_n + \sum_{i=1}^n a_iP_i.
\nn\ee
We claim that $D \in \mathcal{U}(r,t;j_1, \ldots, j_m)$. It follows from the construction that
\be
\begin{array}{rl}
\deg \left(\overline{D}\right) & = \dd (m+1) \left( n-(j_0+ \cdots + j_m)\right) + \sum_{\ell=0}^m (m-\ell) j_\ell \\
&=\dd (m+1)n + \sum_{\ell=0}^m (m-\ell-m-1)j_\ell \\
&=\dd (m+1)n-\sum_{\ell=0}^m (\ell+1) j_\ell =t,
\end{array}
\nn\ee
where we used (\ref{ep2.lemma.non.emptiness.mathcal.u}). Moreover $\deg (D)=\deg \left(\overline{D}\right)+(r-t)=r$, $j_\ell(D)=|S_\ell|=j_\ell$
for each $1 \le \ell \le m$, and hence $D \in \mathcal{U}(r,t;j_1, \ldots, j_m)$.

Next we consider the case $j_m+\cdots+j_1+j_0=n$. This case implies that (cf. (\ref{ep2.lemma.non.emptiness.mathcal.u}))
\be
mn=t+j_1+2j_2+ \cdots + mj_m.
\nn\ee
Therefore we construct $\overline{D}$ similarly and $D$ using the existence of a degree $r-t$ positive divisor
whose support is disjoint
from the set $\{P_1, \ldots, P_n\}$.
\end{proof}

\begin{definition} \label{constant.C.a.b}
For integers $a \ge b \ge 0$ with $b \le n$ and a set $\{Q_1, \ldots, Q_b\}$ of rational places, let $C_{a,b}$ denote the cardinality
of the set of positive divisors given by
\be
\left\{ D \ge 0 : \deg (D)=a, \ \supp \left(\overline{D}\right) = \{Q_1, \ldots, Q_b\} \right\}.
\nn\ee
Note that $C_{a,b}$ is independent of the
choice of the set $\{Q_1, \ldots, Q_b\}$, only the cardinality $b$ of this
set matters.
\end{definition}

\begin{lemma} \label{lemma.cardinality.U}
For $r \ge t\ge 0$, $j_1, \ldots, j_m \ge 0$, and
$mn-(j_1+ \cdots + mj_m) \le t \le (m+1)n- \left(2j_1+ \cdots + (m+1)j_m\right)$, the cardinality of
$\mathcal{U}(r,t;j_1, \ldots, j_m)$ is
\be
\begin{array}{l}
\dd {n \choose j_m} {n-j_m \choose j_{m-1} } \cdots {n-(j_2+j_3+ \cdots + j_m) \choose j_1}
{n-(j_1+j_2+ \cdots +j_m) \choose t-mn+(j_1+2j_2+\cdots+mj_m)}
\\ \\
\dd \times C_{r-mn+(j_1+2j_2+ \cdots + mj_m), t-mn + (j_1+2j_2+ \cdots + mj_m)}.
\end{array}
\nn\ee
\end{lemma}
\begin{proof}
We prove the lemma for $m=2$ and the general case is similar. For $D \in \mathcal{U}(r,t;j_1,j_2)$, let
$S_2=\big\{P \in \{P_1, \ldots, P_n\}: v_P(\overline{D}) =0 \big\}$,
$S_1=\big\{P \in \{P_1, \ldots, P_n\}: v_P(\overline{D})=1\big\}$,
$S_0=\big\{P \in \{P_1, \ldots, P_n\}: v_P(\overline{D})=2\big\}$,
and
$S=\big\{P \in \{P_1, \ldots, P_n\}: v_P(\overline{D}) = 3\big\}=
\big\{P \in \{P_1, \ldots, P_n\}: v_P(D) \ge 3\big\}$.
Note that $|S_2|=j_2$ and $|S_1|=j_1$ and that by (\ref{ep0.lemma.non.emptiness.mathcal.u})
we get $|S_0|=j_0(D)=3n-(2j_1+3j_2)-t$. The choices
of $S_2$, $S_1$, and $S_0$ determine $S$.
We have $|S|=n-(j_1+j_2) - |S_0|=t-2n+(j_1+2j_2)$.
Hence there are
\be
{n \choose j_2} {n-j_2 \choose j_1} {n-(j_1+j_2) \choose t-2n+(j_1+2j_2)}
\nn\ee
choices for these subsets. Assume that the subsets $S_2$, $S_1$, $S_0$, and $S$ are determined.
For a corresponding $D \in \mathcal{U}(r,t;j_1,j_2)$, let $D_1=b_1P_1+ \cdots +b_n P_n$, where
\be
b_i=\left\{\begin{array}{ll}
v_{P_i}(D)=v_{P_i}(\overline{D}) & \mbox{if $P_i \in S_2 \cup S_1 \cup S_0$}, \\
2=v_{P_i}(\overline{D})-1 & \mbox{if $P_i \in S$}.
\end{array}
\right.
\nn\ee
Moreover let $E=D-D_1$. Then $E$ is a positive divisor and $\supp \left(\overline{E}\right)=S$.
Note that
\be
\deg \left( D_1 \right)=t-|S|, \;\;  \deg (E)=\deg (D)-\deg \left(D_1\right)=r-t+|S|.
\nn\ee
Hence
\be
|\supp \left(\overline{E} \right)| = t-2n + (j_1+2j_2), \;\;
\deg (E)=r-2n+(j_1+2j_2).
\nn\ee
Using Definition \ref{constant.C.a.b}, we obtain that there are $C_{r-2n+(j_1+2j_2),t-2n+(j_1+2j_2)}$
choices for $E$, which completes the proof.
\end{proof}

Recall that for integers $r\ge s \ge 0$ and nonnegative integers $X_1, \ldots, X_m$, the set
$\mathcal{V}_m(r,s;X_1, \ldots, X_m)$ is defined in Definition
\ref{definition.V.m}. Using  Definition
\ref{definition.set.U} and Lemma \ref{lemma.non.emptiness.mathcal.u},
we can write the set $\mathcal{V}_m(r,s;X_1, \ldots, X_m)$ as the disjoint union
\be \label{disjoint.V.m}
\dd \mathcal{V}_m(r,s;X_1, \ldots, X_m)= \bigsqcup_{j_m}
\bigsqcup_{j_{m-1}}
\cdots
\bigsqcup_{j_1}
\bigsqcup_{t} \mathcal{U}(r,t;j_1, \ldots,j_m),
\ee
where the $m$-tuples $(j_1, \ldots, j_m)$ of indices run over the finite set of $m$-tuples of integers
satisfying
\be \label{running.space.j.m.tuple}
\begin{array}{c}
\begin{array}{l}
0 \le j_m \le 2X_m, \; 0 \le j_{m-1} \le 2X_{m-1}+X_m, \ldots, \\
0 \le j_1 \le 2X_1 + \sum_{\nu=2}^m X_\nu,
\end{array} \\
2j_1+3j_2+ \cdots+(m+1)j_m \le 2 (2X_1+3X_2+ \cdots+ (m+1)X_m),
\end{array}
\ee
and for each $m$-tuple satisfying (\ref{running.space.j.m.tuple}), the index $t$  runs from
$\max\left(s,mn-(j_1+2j_2+ \cdots +mj_m)\right)$ to $\min\left(r,(m+1)n - (2j_1+3j_2+ \cdots + (m+1)j_m)\right)$.

Combining (\ref{disjoint.V.m}) and Lemma \ref{lemma.cardinality.U}, we can compute the cardinality of
the set $\mathcal{V}_m(r,s;X_1,\ldots,X_m)$.

\section{Asymptotic Upper Bound on the Cardinality of $\mathcal{V}_1(r,s;X_1)$ } \label{section.asymptotic.size.V.m.case.m.1}

In this section we obtain an asymptotic upper bound on the cardinality of $\mathcal{V}_m(r,s;X_1, \ldots, X_m)$ for the case $m=1$
in a suitable sequence of global function fields (see Corollary \ref{corollary.asymptotic.card.V.1}). The assumption $m=1$ is made for simplicity and for the clarity of the exposition. Later in Section \ref{section.asymptotic.size.V.m.case.m.general}
we generalize this asymptotic upper bound to the case $m \ge 1$.

The asymptotic upper bound for the cardinality of $\mathcal{V}_m(r,s;X_1, \ldots, X_m)$ will be used later
to prove the existence of a sequence of distinguished divisors on the basis of Proposition \ref{proposition.existence.G}.

\begin{definition} \label{definition.E}
Let $E$ be the real-valued function defined on the interval $[0,1]$ as follows: for $0<x<1$ we put
$E(x)=-x\log_q x- (1-x) \log_q (1-x)$ and for $x \in \{0,1\}$ we put $E(0)=E(1)=\lim_{x\ra 0^+} E(x)= \lim_{x \ra 1^-} E(x)=0$.
\end{definition}

Using Stirling's formula, we obtain the following well-known results. For any real number $0 \le \alpha \le 1$, we have
\be \label{asymptotic.binomial.alpha}
\lim_{n \ra \infty} \frac{\log_q {n \choose \lfloor \alpha n \rfloor}}{n}=E(\alpha).
\ee
For any real numbers  $0 \le \alpha_1 \le 1$ and $0 \le \alpha_2 < 1$ with $\alpha_1+\alpha_2 \le 1$, we have
\be \label{asymptotic.binomial.alpha1.alpha2}
\lim_{n \ra \infty} \frac{\log_q {n - \lfloor \alpha_2 n \rfloor \choose \lfloor \alpha_1 n \rfloor}}{n}
=(1-\alpha_2)E\left( \frac{\alpha_1}{1-\alpha_2}\right).
\ee

Now we state an important assumption and introduce related notation.

\begin{description}

 \item[{\bf Assumption 1}] Assume that $\left({F_i/\F_q}\right)_{i=1}^\infty$ is a sequence of global function fields with full constant field $\F_q$, with $g_i \ra \infty$
as $i \ra \infty$, and with $\lim_{i \ra \infty} \frac{n_i}{g_i}=\gamma >0$, where
$n_i$ and $g_i$ denote the number of rational places and the genus of $F_i$,
respectively.

\end{description}

We will use the following proposition in our upper bounds.

\begin{proposition} \label{proposition.C.a.b}
Under Assumption 1,
let $(a_i)_{i=1}^\infty$ and $(b_i)_{i=1}^\infty$ be sequences of integers such that $a_i  \ge b_i \ge 0$ and $b_i \le n_i$
for all $i \ge 1$.
We also assume that there exist the limits
\be \label{Assumption2}
\lim_{i \ra \infty} \frac{a_i}{n_i}=a, \ \ \lim_{i \ra \infty} \frac{b_i}{n_i}=b \ \ \mbox{with }
0<b\le a < \infty.
\ee
For each $i \ge 1$, let $C_{a_i,b_i}^{(i)}$ denote the cardinality of the set of positive
divisors given in Definition \ref{constant.C.a.b}
for a suitable set $\{Q^{(i)}_1, \ldots, Q^{(i)}_{b_i}\}$ of rational places of $F_i$. Then we have
\be
\limsup_{i \ra \infty} \frac{\log_q C_{a_i,b_i}^{(i)}}{n_i} \le
\left\{
\begin{array}{ll}
a E\left(\frac{b}{a}\right) & \mbox{if $\frac{b}{a} \ge 1-\frac{1}{q}$}, \\
a-b \log_q (q-1) & \mbox{if $\frac{b}{a} \le 1-\frac{1}{q}$}.
\end{array}
\right.
\nn\ee
\end{proposition}
\begin{proof}
This follows from Definition \ref{constant.C.a.b} and the proof of \cite[Lemma 3.4.10]{TV}.
\end{proof}

Let $y,\sigma,x_1 \ge 0$ be real numbers. Under Assumption 1, for each $i \ge 1$  we define the integers
\be \label{definition.r.i.s.i.W.i}
r_i=\left\lfloor \left(1+y+\frac{\sigma}{\gamma}\right) n_i \right\rfloor, \; s_i=\lfloor (1+y)n_i \rfloor, \; X_1^{(i)}=\lfloor x_1 n_i \rfloor.
\ee

Let $\mathcal{V}_1^{(i)}(r_i,s_i;X_1^{(i)})$ be the set of positive divisors of degree $r_i$ of $F_i$, which is defined using Definition
\ref{definition.V.m}. We note that for each real number $0 \le t_1 \le 2x_1$ and each integer $i\ge 1$, we have
\be
\max\{s_i, n_i - \lfloor t_1 n_i \rfloor\}=s_i.
\nn\ee
Moreover, if
\be \label{Assumption3}
1+y+\frac{\sigma}{\gamma} < 2 -4x_1 \ \ \mbox{or equivalently} \ \ y + 4x_1 + \frac{\sigma}{\gamma} < 1
\ee
holds, then for each real number $0 \le t_1 \le 2x_1$ and integer $i\ge 1$ we also have
\be
\min\{r_i,2n_i-2\lfloor t_1 n_i \rfloor \}=r_i.
\nn\ee

\begin{definition} \label{definition.S}
For real numbers $y>0$, $x_1,\sigma\ge 0$ satisfying (\ref{Assumption3}) and real numbers  $0 \le t_1 \le 2x_1$,
$0 \le x \le \frac{\sigma}{\gamma}$, let $S(\sigma,y,x,t_1)$
be the real-valued function
\be
\begin{array}{rl}
\dd S(\sigma,y,x,t_1)= & \dd E(t_1) + (1-t_1)E\left(\frac{y+x+t_1}{1-t_1}\right) \\ \\
& + \left\{
\begin{array}{ll}
\left(y+ \frac{\sigma}{\gamma} + t_1\right) E\left( \frac{y+x+t_1}{y+ \frac{\sigma}{\gamma} + t_1}\right) &
\mbox{if $ \frac{y+x+t_1}{y+\frac{\sigma}{\gamma} + t_1} \ge 1 -\frac{1}{q}$}, \\ \\
\left(y+\frac{\sigma}{\gamma}+t_1\right) - (y+x+t_1)\log_q (q-1) &
\mbox{if $ \frac{y+x+t_1}{y+\frac{\sigma}{\gamma} + t_1} \le 1 -\frac{1}{q}$}.
\end{array}
\right.
\end{array}
\nn\ee
Note that by (\ref{Assumption3})  we have $4x_1 <1$ and hence $t_1 < \frac{1}{2}$.
\end{definition}

\begin{proposition} \label{proposition.asymtotic.cardinality.U.1}
Under Assumption 1, let $y>0$ and $x_1,\sigma \ge 0$ be real numbers satisfying (\ref{Assumption3}). For each integer $i \ge 1$
and real numbers $0 \le t_1 \le 2x_1$, $0 \le x \le \frac{\sigma}{\gamma}$, let
$\mathcal{U}^{(i)}\left(\left\lfloor (1+y+\frac{\sigma}{\gamma})n_i\right\rfloor,
\left\lfloor (1+y+x)n_i\right\rfloor; \left\lfloor t_1 n_i \right\rfloor \right)$ be the set of positive divisors of $F_i$
defined in Definition \ref{definition.set.U} for $m=1$. Then for the cardinalities of these sets we have
\be
\limsup_{i \ra \infty}
\frac{\log_q \left| \mathcal{U}^{(i)}\left(\left\lfloor (1+y+\frac{\sigma}{\gamma})n_i\right\rfloor,
\left\lfloor (1+y+x)n_i\right\rfloor; \left\lfloor t_1 n_i \right\rfloor \right) \right|}{n_i} \le S(\sigma,y,x,t_1).
\nn\ee
\end{proposition}
\begin{proof}
Note that $n_i - \lfloor t_1 n_i \rfloor \le \lfloor (1+y+x) n_i \rfloor$ and using (\ref{Assumption3}) we get
$\lfloor (1+y+x)n_i\rfloor \le 2n_i - 2 \lfloor t_1 n_i \rfloor$ for each $x$ and $t_1$ in the range under consideration.
Hence using Lemma \ref{lemma.cardinality.U}, we obtain
\be \label{ep1.proposition.asymtotic.cardinality.U.1}
\begin{array}{ll}
&  \left|\mathcal{U}^{(i)}\left(\left\lfloor \left(1+y+\frac{\sigma}{\gamma}\right)n_i\right\rfloor,
\left\lfloor (1+y+x)n_i\right\rfloor; \left\lfloor t_1 n_i \right\rfloor \right)\right| \\ \\
= & {n \choose \lfloor t_1 n_i \rfloor} { n- \lfloor t_1 n_i \rfloor \choose \lfloor (1+y+x) n_i \rfloor -n_i + \lfloor t_1 n_i \rfloor} \\ \\
&  \times
C^{(i)}_{\lfloor (1+y+\frac{\sigma}{\gamma})n_i \rfloor-n_i+\lfloor t_1 n_i \rfloor, \lfloor (1+y+x)n_i \rfloor - n_i + \lfloor t_1 n_i \rfloor}.
\end{array}
\ee
Using (\ref{asymptotic.binomial.alpha}) and (\ref{asymptotic.binomial.alpha1.alpha2}), we obtain
\be \label{ep2.proposition.asymtotic.cardinality.U.1}
\begin{array}{l}
\dd \lim_{i \ra \infty} \frac{ \log_q {n_i \choose \lfloor t_1 n_i \rfloor}}{n_i} = E(t_1), \\ \\
\dd \lim_{i \ra \infty} \frac{\log_q {n_i - \lfloor t_1 n_i \rfloor \choose \lfloor (1+y+x)n_i \rfloor - n_i + \lfloor t_1 n_i \rfloor}}{n_i}
= (1-t_1)E\left(\frac{y+x+t_1}{1-t_1}\right).
\end{array}
\ee
Note that $\lim_{i \ra \infty} \frac{\lfloor(1+y+\frac{\sigma}{\gamma})n_i \rfloor -n_i + \lfloor t_1 n_i \rfloor}{n_i}=y+\frac{\sigma}{\gamma}+t_1$
and $\lim_{i \ra \infty} \frac{\lfloor (1+y+x)n_i\rfloor -n_i + \lfloor t_1 n_i \rfloor}{n_i}=y+x+t_1$. Hence from
Proposition \ref{proposition.C.a.b} we get
\be \label{ep3.proposition.asymtotic.cardinality.U.1}
\begin{array}{l}
\dd \limsup_{i \ra \infty} \frac{\log_q C^{(i)}_{\lfloor (1+y+\frac{\sigma}{\gamma})n_i \rfloor -n_i + \lfloor t_1 n_i \rfloor, \lfloor (1+y+x) n_i \rfloor -n_i + \lfloor t_1 n_i \rfloor}}{n_i} \\ \\
\le
\dd
\left\{
\begin{array}{ll}
(y+\frac{\sigma}{\gamma} + t_1) E\left(\frac{y+x+t_1}{y+\frac{\sigma}{\gamma}+t_1}\right) & \mbox{if $\frac{y+x+t_1}{y+\frac{\sigma}{\gamma}+t_1} \ge 1 -\frac{1}{q}$}, \\
(y+\frac{\sigma}{\gamma}+t_1) - (y+x+t_1)\log_q (q-1) & \mbox{if $\frac{y+x+t_1}{y+\frac{\sigma}{\gamma}+t_1} \le 1 -\frac{1}{q}$.}
\end{array}
\right.
\end{array}
\ee
Using (\ref{ep1.proposition.asymtotic.cardinality.U.1}), (\ref{ep2.proposition.asymtotic.cardinality.U.1}), (\ref{ep3.proposition.asymtotic.cardinality.U.1}),
and Definition \ref{definition.S}, we complete the proof.
\end{proof}

\begin{corollary} \label{corollary.asymptotic.card.V.1}
Under Assumption 1, let $y>0$ and $x_1,\sigma \ge 0$ be real numbers satisfying (\ref{Assumption3}). For each integer $i \ge 1$,
let $r_i, s_i$, and $X_1^{(i)}$ be the integers defined in (\ref{definition.r.i.s.i.W.i}) and let
$\mathcal{V}_1^{(i)}(r_i,s_i;X_1^{(i)})$ be the set of positive divisors of $F_i$ defined in Definition \ref{definition.V.m} for $m=1$.
Then for the cardinalities of these sets we have
\be
\limsup_{i \ra \infty} \frac{\log_q |\mathcal{V}_1^{(i)}(r_i,s_i;X_1^{(i)})|}{n_i} \le \max S(\sigma,y,x,t_1),
\nn\ee
where the maximum is over all real numbers $x$ and $t_1$ satisfying $0 \le x \le \frac{\sigma}{\gamma}$ and $0 \le t_1 \le 2x_1$.
\end{corollary}
\begin{proof}
Using  (\ref{disjoint.V.m}) and Lemma \ref{lemma.cardinality.U} for each $i \ge 1$, we obtain that
\be \label{ep1.corollary.asymptotic.card.V.1}
|\mathcal{V}_1^{(i)}(r_i,s_i;X_1^{(i)})|=\sum_{j_1=0}^{ 2X_1^{(i)} } \sum_{t} |\mathcal{U}^{(i)}(r_i,t;j_1)|,
\ee
where $t$ runs from $\max\{s_i,n_i-j_1\}$ to $\min\{r_i,2n_i-2j_1\}$. Note that $s_i \ge n_i-j_1$
and $r_i \le 2n_i - 2j_1$ for each $i \ge 1$ and $0 \le j_1 \le 2X_1^{(i)}$. Moreover, for the number of
terms $\left(2X_1^{(i)} +1\right)\left( r_i-s_i+1\right)$ in the summation in (\ref{ep1.corollary.asymptotic.card.V.1}) we have
\be
\begin{array}{l}
\dd \lim_{i \ra \infty} \frac{\log_q \left(\left( 2X_1^{(i)} +1\right) \left(r_i-s_i+1\right)\right)}{n_i} \\ \\
\dd =\lim_{i \ra \infty} \left\{
\frac{\log_q \left(2x_1 + 1/n_i\right) + \log_q\left(\frac{\sigma}{\gamma} +1/n_i\right)}{n_i} + 2 \frac{\log_q n_i}{n_i}\right\}=0.
\end{array}
\nn\ee
Therefore, using the method of the proof of \cite[Proposition 4.3]{NO2} and Proposition \ref{proposition.asymtotic.cardinality.U.1},
we complete the proof.
\end{proof}

\begin{definition} \label{definition.I.w}
Under Assumption 1, let $y >0$ and $x_1 \ge 0$ be real numbers such that $y+4x_1 < 1$. For  $\sigma \ge 0$
and $y+4x_1+ \frac{\sigma}{\gamma} < 1$, let $I_{y,x_1}(\sigma)$ be the real-valued function of $\sigma$ defined by
\be
I_{y,x_1}(\sigma)=\max S(\sigma,y,x,t_1),
\nn\ee
where the maximum is over all real numbers $x$ and $t_1$ such that $0 \le t_1 \le 2x_1$ and $0 \le x \le \frac{\sigma}{\gamma}$.
\end{definition}

By straightforward manipulations, the expression for $S(\sigma,y,x,t_1)$ is simplified to
{\footnotesize
\be \label{simplified.S.w}
\begin{array}{l}
\dd S(\sigma,y,x,t_1) \\ \\
\dd =-t_1\log_q t_1 \\ \\
\dd -(y+x+t_1)\log_q(y+x+t_1) \\ \\
\dd -(1-y-x-2t_1)\log_q(1-y-x-2t_1) \\ \\
\dd +\left\{
\begin{array}{l}
\dd -(y+x+t_1)\log_q(y+x+t_1) - \left(\frac{\sigma}{\gamma}-x\right) \log_q\left(\frac{\sigma}{\gamma}-x\right) \\ \\
\dd + \left(y+ \frac{\sigma}{\gamma} + t_1\right) \log_q \left(y+ \frac{\sigma}{\gamma} + t_1\right)
 \hspace{2cm} \mbox{if} \;\; \frac{y+x+t_1}{y+\frac{\sigma}{\gamma}+t_1} \ge 1 -\frac{1}{q},  \\ \\ \\
\dd \left(y+ \frac{\sigma}{\gamma} + t_1\right) - (y+ x + t_1) \log_q (q-1)
 \hspace{1.4cm} \mbox{if} \;\; \frac{y+x+t_1}{y+\frac{\sigma}{\gamma}+t_1} \le 1 -\frac{1}{q}.
\end{array} \right.
\end{array}
\ee
}

We first show that $I_{y,x_1}(\sigma)$ is a strictly increasing function of $\sigma$.

\begin{lemma} \label{lemma.I.is.increasing.of.sigma}
Under the assumptions of Definition \ref{definition.I.w},
the real-valued function $I_{y,x_1}(\sigma)$ is a strictly
increasing function of $\sigma$ on its domain of definition, which is the interval of $\sigma$ such that
$\sigma \ge 0$ and $y+4x_1+\frac{\sigma}{\gamma} < 1$.
\end{lemma}
\begin{proof}
Using the expression (\ref{simplified.S.w}), for the
partial derivative of $S(\sigma,y,x,t_1)$ with respect to $\sigma$ we obtain
\be
\frac{\partial S}{\partial \sigma}(\sigma,y,x,t_1)
=\left\{
\begin{array}{ll}
\frac{1}{\gamma} \log_q \frac{y + \frac{\sigma}{\gamma}+t_1}{\frac{\sigma}{\gamma}-x} &
\mbox{if $\frac{y+x+t_1}{y+\frac{\sigma}{\gamma}+t_1} \ge 1-\frac{1}{q}$}, \\
\frac{1}{\gamma} & \mbox{if $\frac{y+x+t_1}{y+\frac{\sigma}{\gamma}+t_1} \le 1-\frac{1}{q}$}.
\end{array}
\right.
\nn\ee
Therefore  $\frac{\partial S}{\partial \sigma}(\sigma,y,x,t_1) > 0$ for each $0 \le x < \frac{\sigma}{\gamma}$
and $0 \le t_1 \le 2x_1$. Moreover $\lim_{x \ra \frac{\sigma}{\gamma}^-} \frac{\partial S}{\partial \sigma}(\sigma,y,x,t_1)=+\infty$
for $0 \le t_1 \le 2x_1$. This completes the proof.
\end{proof}

\begin{lemma} \label{lemma.del.S.w.to.t1.x}
Under the assumptions of Definition \ref{definition.I.w}, for the partial
derivatives $\frac{\partial S}{\partial t_1}(\sigma,y,x,t_1)$ and $\frac{ \partial S}{\partial x}(\sigma, y,x,t_1)  $
of $S(\sigma,y,x,t_1)$ with respect to $t_1$ and $x$ we obtain
\be
 \frac{\partial S}{\partial t_1}(\sigma,y,x,t_1)
=
\log_q \frac{(1-y-x-2t_1)^2}{t_1(y+x+t_1)}
+
\left\{
\begin{array}{ll}
\log_q \frac{y+\frac{\sigma}{\gamma} + t_1}{y+x+t_1} & \mbox{if $\frac{y+x+t_1}{y+\frac{\sigma}{\gamma}+t_1} \ge 1 - \frac{1}{q}$}, \\
\log_q\frac{q}{q-1} & \mbox{if $\frac{y+x+t_1}{y+\frac{\sigma}{\gamma}+t_1} \le 1 - \frac{1}{q}$},
\end{array}
\right.
\nn\ee
and
\be
 \frac{\partial S}{\partial x}(\sigma,y,x,t_1)
=
\log_q \frac{1-y-x-2t_1}{y+x+t_1}
+
\left\{
\begin{array}{ll}
\log_q \frac{\frac{\sigma}{\gamma}-x}{y+x+t_1} & \mbox{if $\frac{y+x+t_1}{y+\frac{\sigma}{\gamma}+t_1} \ge 1 - \frac{1}{q}$}, \\
-\log_q (q-1) & \mbox{if $\frac{y+x+t_1}{y+\frac{\sigma}{\gamma}+t_1} \le 1 - \frac{1}{q}$}.
\end{array}
\right.
\nn\ee
\end{lemma}
\begin{proof}
Let $S_1$, $T_1$, and $T_2$ denote the following expressions from (\ref{simplified.S.w}):
\be
\begin{array}{rl}
\dd S_1= & -t_1\log_q t_1 - (y+x+t_1) \log_q(y+x+t_1) \\
&- (1-y-x-2t_1)\log_q(1-y-x-2t_1), \\ \\
\dd T_1=& -(y+x+t_1) \log_q (y+x+t_1) - \left(\frac{\sigma}{\gamma} -x\right) \log_q \left(\frac{\sigma}{\gamma}-x\right) \\
&+ \left(y+ \frac{\sigma}{\gamma} + t_1\right)\log_q(y+\frac{\sigma}{\gamma}+t_1), \\ \\
\dd T_2= & \left(y+\frac{\sigma}{\gamma}+t_1\right) - (y+x+t_1) \log_q (q-1).
\end{array}
\nn\ee
For their partial derivatives with respect to $t_1$ and $x$ we obtain
\be
\begin{array}{l}
\dd \frac{\partial S_1}{\partial t_1}= -\log_q t_1 - \log_q (y+x+t_1) + 2 \log_q (1-y-x-2t_1), \\ \\
\dd \frac{\partial T_1}{\partial t_1}=-\log_q (y+x+t_1) + \log_q\left(y+ \frac{\sigma}{\gamma}+t_1\right), \\ \\
\dd \frac{\partial T_2}{\partial t_1}=1-\log_q(q-1)=\log_q\frac{q}{q-1},
\end{array}
\nn\ee
and
\be
\begin{array}{l}
\dd \frac{\partial S_1}{\partial x}= - \log_q(y+x+t_1) + \log_q(1-y-x-2t_1), \\ \\
\dd \frac{\partial T_1}{\partial x}=-\log_q(y+x+t_1) + \log_q\left(\frac{\sigma}{\gamma}-x\right), \\ \\
\dd \frac{\partial T_2}{\partial x}=-\log_q(q-1).
\end{array}
\nn\ee
Using (\ref{simplified.S.w}) and combining the partial derivatives above, we
get the desired formulas.
\end{proof}

\begin{corollary} \label{corollary.I.in.restricted.case}
Under the assumptions of Definition \ref{definition.I.w}, furthermore if all of the following
conditions
\begin{enumerate}
\item[C1:] $\frac{\sigma}{\gamma} \le \frac{y}{q-1}$,
\item[C2:] $2x_1 \left(y+\frac{\sigma}{\gamma}+2x_1\right)^2
< \left( 1-y -\frac{\sigma}{\gamma}-4x_1\right)^2 \left( y + \frac{\sigma}{\gamma}\right)$,
\item[C3:] $\frac{\sigma}{\gamma} (1-y) < y^2$,
\end{enumerate}
hold, then we have
\be
\begin{array}{l}
I_{y,x_1}(\sigma)=S\left(\sigma,y,0,2x_1\right) \\ \\
=E(2x_1) + \left(1-2x_1\right) E\left( \frac{y+2x_1}{1-2x_1}\right) +
\left(y+ \frac{\sigma}{\gamma}+2x_1\right)E\left( \frac{y+2x_1}{y+\frac{\sigma}{\gamma}+2x_1}\right).
\end{array}
\nn\ee
\end{corollary}
\begin{proof}
Assume that $0 \le x \le \frac{\sigma}{\gamma}$ and $0 \le t_1 \le 2x_1$. First we observe that
\be
\frac{y+x+t_1}{y+\frac{\sigma}{\gamma}+t_1} \ge \frac{y}{y+\frac{\sigma}{\gamma}}.
\nn\ee
Using condition C1 we obtain
\be \label{ep1.corollary.I.in.restricted.case}
\frac{y+x+t_1}{y+\frac{\sigma}{\gamma}+t_1} \ge \frac{y}{y+\frac{\sigma}{\gamma}} \ge 1-\frac{1}{q}.
\ee
Moreover, using condition C2 we also get
\be
\begin{array}{rl}
t_1(y+x+t_1)^2 & \le 2x_1 \left(y+\frac{\sigma}{\gamma} + 2x_1 \right)^2
 < \left(y+ \frac{\sigma}{\gamma}\right) \left(1-y-\frac{\sigma}{\gamma}-4x_1 \right)^2 \\ \\
& \le \left(y+\frac{\sigma}{\gamma}+t_1\right) \left(1-y-x-2t_1\right)^2.
\end{array}
\nn\ee
Therefore by Lemma \ref{lemma.del.S.w.to.t1.x} and (\ref{ep1.corollary.I.in.restricted.case}) we have
$\frac{\partial S}{\partial t_1}(\sigma,y,x,t_1) > 0$. Similarly, condition C3 implies
\be
\left(\frac{\sigma}{\gamma}-x\right)(1-y-x-2t_1) \le \frac{\sigma}{\gamma}(1-y) < y^2 \le (x+y+t_1)^2,
\nn\ee
and by Lemma \ref{lemma.del.S.w.to.t1.x} we also have $\frac{\partial S}{\partial x}(\sigma,y,x,t_1) < 0$.
Hence we obtain $I_{y,x_1}(\sigma)=S\left(\sigma,y,0,2x_1\right)$. Using Definition \ref{definition.S} we complete
the proof.
\end{proof}

\section{Asymptotic Upper Bound on the Cardinality of $\mathcal{V}_m(r,s;X_1, \ldots, X_m)$ for the General Case $m \ge 1$} \label{section.asymptotic.size.V.m.case.m.general}

In this section we obtain generalizations of the results of Section \ref{section.asymptotic.size.V.m.case.m.1} to the general case
$m \ge 1$. For simplicity we begin with the case $m=2$, which corresponds to the two-variable case $t_1, t_2$.

\begin{definition} \label{definition.S.case.m.2}
Let $\gamma > 0$ be as in Assumption 1 (cf. Section \ref{section.asymptotic.size.V.m.case.m.1}).
Let $y>0$, $x_1,x_2,\sigma \ge 0$ be real numbers satisfying
\be \label{Assume.condition.on.y.x1.x2.sigma.m.2}
y+ 2(2x_1+3x_2) + \frac{\sigma}{\gamma} < 1.
\ee
For real numbers  $0 \le x \le \frac{\sigma}{\gamma}$ and $0 \le t_1, t_2$ satisfying $t_2 \le 2 x_2$, $t_1 \le 2x_1 + x_2$,
and $2t_1 + 3t_2 \le 2(2x_1+3x_2)$,
let $S(\sigma,y,x,t_1,t_2)$
be the real-valued function
\be
\begin{array}{l}
\dd S(\sigma,y,x,t_1,t_2)=  E(t_2) + (1-t_2) E \left( \frac{t_1}{1-t_2} \right) \\ \\
 \dd + (1-t_1-t_2)E\left(\frac{y+x+t_1+2t_2}{1-t_1-t_2}\right) \\ \\
 + \left\{
\begin{array}{ll}
\left(y+ \frac{\sigma}{\gamma} + t_1+2t_2\right) E\left( \frac{y+x+t_1+2t_2}{y+ \frac{\sigma}{\gamma} + t_1+2t_2}\right) &
\mbox{if $ \frac{y+x+t_1+2t_2}{y+\frac{\sigma}{\gamma} + t_1+2t_2} \ge 1 -\frac{1}{q}$}, \\ \\
\left(y+\frac{\sigma}{\gamma}+t_1+2t_2\right) - (y+x+t_1+2t_2)\log_q (q-1) &
\mbox{if $ \frac{y+x+t_1+2t_2}{y+\frac{\sigma}{\gamma} + t_1+2t_2} \le 1 -\frac{1}{q}$}.
\end{array}
\right.
\end{array}
\nn\ee
Note that by (\ref{Assume.condition.on.y.x1.x2.sigma.m.2}) we have $2(2x_1+3x_2) <1$ and hence $t_1 +t_2  \le t_1 + \frac{3}{2} t_2 < \frac{1}{2}$.
\end{definition}

Instead of stating a generalization of Proposition  \ref{proposition.asymtotic.cardinality.U.1} explicitly, we
prefer to give a generalization of Corollary \ref{corollary.asymptotic.card.V.1} directly in the following proposition,
whose proof includes a generalization of Proposition \ref{proposition.asymtotic.cardinality.U.1}.

\begin{proposition} \label{propostion.asymptotic.card.V.m.case.2}
Under Assumption 1 (cf. Section \ref{section.asymptotic.size.V.m.case.m.1}), let
$y > 0$ and $x_1,x_2,\sigma \ge 0$ be real numbers satisfying (\ref{Assume.condition.on.y.x1.x2.sigma.m.2}).
For each integer $i \ge 1$,
let $r_i=\left\lfloor \left(2 + y + \frac{\sigma}{\gamma} \right) n_i \right \rfloor$,
$s_i=\left\lfloor (2+y)n_i \right \rfloor$, $X_1^{(i)}=\left \lfloor x_1n_i \right \rfloor$,
$X_2^{(i)}=\left \lfloor x_2 n_i \right\rfloor$, and
$\mathcal{V}^{(i)}_2(r_i,s_i;X_1^{(i)},X_2^{(i)})$ be the set of positive divisors of $F_i$
defined in Definition \ref{definition.V.m} for
$m=2$. Then for the cardinalities of these sets we have
\be
\limsup_{i \ra \infty} \frac{ \log_q \left| \mathcal{V}_2^{(i)}(r_i,s_i;X_1^{(i)},X_2^{(i)})\right|}{n_i}
\le \max S(\sigma,y,x,t_1,t_2),
\nn\ee
where the maximum is over all real numbers $x$ and $t_1, t_2$ satisfying $0 \le x \le \frac{\sigma}{\gamma}$
and $0 \le t_2 \le 2 x_2$,  $0 \le t_1 \le 2x_1+x_2$, and  $2t_1 + 3t_2 \le 2(2x_1+3x_2)$.
\end{proposition}
\begin{proof}
We follow similar methods as in the proofs of Proposition \ref{proposition.asymtotic.cardinality.U.1} and
Corollary \ref{corollary.asymptotic.card.V.1}. First note that for each integer $i \ge 1$ and real numbers
$0 \le t_1, t_2$ with $2t_1+3t_2 \le 2(2x_1+3x_2)$, using (\ref{Assume.condition.on.y.x1.x2.sigma.m.2}) we obtain
$r_i \le 3n_i - \left(2 \lfloor t_1 n_i\rfloor + 3 \lfloor t_2 n_i \rfloor \right)$. Moreover it is also clear that
$s_i \ge 2n_i -\left( \lfloor t_1 n_i \rfloor + 2 \lfloor t_2 n_i \rfloor \right)$ for each integer $i \ge 1$ and
real numbers $t_1, t_2 \ge 0$. Hence using (\ref{disjoint.V.m}) and Lemma \ref{lemma.cardinality.U} as in the proof of
Corollary \ref{corollary.asymptotic.card.V.1}, for integers $i \ge 1$ and
real numbers $0 \le x, t_1, t_2$ such that $x \le \frac{\sigma}{\gamma}$,
$t_2 \le 2 x_2$, $t_1 \le 2x_1+x_2$, and  $2t_1 + 3t_2 \le 2(2x_1+3x_2)$,
we need to consider the cardinality
$ \left|\mathcal{U}^{(i)}\left(r_i, \lfloor (2+y+x) n_i \rfloor; \lfloor t_1 n_i \rfloor, \lfloor t_2 n_i \rfloor \right) \right|$
of the set of positive divisors of $F_i$ defined in Definition \ref{definition.set.U} for $m=2$. By Lemma \ref{lemma.cardinality.U} we have
\be
\begin{array}{ll}
&  \left|\mathcal{U}^{(i)}\left(r_i,
\left\lfloor (2+y+x)n_i\right\rfloor; \left\lfloor t_1 n_i \right\rfloor, \left\lfloor t_2 n_i \right \rfloor
\right)\right| \\ \\
= & {n \choose \lfloor t_2 n_i \rfloor} {n - \lfloor t_2 n_i \rfloor \choose \ \lfloor t_1 n_i \rfloor}
{ n- \left(\lfloor t_1 n_i \rfloor + \lfloor t_2 n_i \rfloor \right) \choose
\lfloor (2+y+x) n_i \rfloor -2n_i + \left( \lfloor t_1 n_i \rfloor + 2 \lfloor t_2 n_i \rfloor \right)} \\ \\
&  \times
C^{(i)}_{r_i-2n_i+\lfloor t_1 n_i \rfloor + 2 \lfloor t_2 n_i \rfloor,
\lfloor (2+y+x)n_i \rfloor - 2n_i + \lfloor t_1 n_i \rfloor + 2 \lfloor t_2 n_i \rfloor}.
\end{array}
\nn\ee
We complete the proof using similar arguments as in the proofs of Proposition \ref{proposition.asymtotic.cardinality.U.1}
and Corollary \ref{corollary.asymptotic.card.V.1}.
\end{proof}

Now we generalize Definition \ref{definition.I.w}.

\begin{definition} \label{definition.I.w.case.m.2}
Under Assumption 1 (cf. Section \ref{section.asymptotic.size.V.m.case.m.1}), let $y>0$ and $x_1,x_2 \ge 0$ be real numbers
such that $y+2(2x_1+3x_2) < 1$. For $\sigma \ge 0$ and $y+2(2x_1+3x_2)+\frac{\sigma}{\gamma} <1$, let
$I_{y,x_1,x_2}(\sigma)$ be the real-valued function of $\sigma$ defined by
\be
I_{y,x_1,x_2}(\sigma)=\max S(\sigma,y,x,t_1,t_2),
\nn\ee
where the maximum is over all real numbers $x, t_1$, and $t_2$ such that $0 \le x \le \frac{\sigma}{\gamma}$ and
 $0 \le t_2 \le 2 x_2$,  $0 \le t_1 \le 2x_1+x_2$, and  $2t_1 + 3t_2 \le 2(2x_1+3x_2)$.
\end{definition}

The following lemma generalizes Lemma \ref{lemma.I.is.increasing.of.sigma}.

\begin{lemma} \label{lemma.I.increasing.case.m.2}
Under the assumptions of Definition \ref{definition.I.w.case.m.2},
the real-valued function $I_{y,x_1,x_2}(\sigma)$ is a strictly
increasing function of $\sigma$ on its domain of definition, which is the interval of $\sigma$ such that
$\sigma \ge 0$ and $y+2(2x_1+3x_2)+\frac{\sigma}{\gamma} < 1$.
\end{lemma}
\begin{proof}
For the
partial derivative of $S(\sigma,y,x,t_1,t_2)$ with respect to $\sigma$ we obtain
\be
\frac{\partial S}{\partial \sigma}(\sigma,y,x,t_1,t_2)
=\left\{
\begin{array}{ll}
\frac{1}{\gamma} \log_q \left( \frac{y + \frac{\sigma}{\gamma}+t_1+2t_2}{\frac{\sigma}{\gamma}-x}\right) &
\mbox{if $\frac{y+x+t_1+2t_2}{y+\frac{\sigma}{\gamma}+t_1+2t_2} \ge 1-\frac{1}{q}$}, \\
\frac{1}{\gamma} & \mbox{if $\frac{y+x+t_1+2t_2}{y+\frac{\sigma}{\gamma}+t_1+2t_2} \le 1-\frac{1}{q}$}.
\end{array}
\right.
\nn\ee
Then the proof is similar to the proof of Lemma \ref{lemma.I.is.increasing.of.sigma}.
\end{proof}

Now we give a generalization of Corollary \ref{corollary.I.in.restricted.case} in the following proposition.

\begin{proposition} \label{proposition.I.in.restricted.case.case.m.2}
Under the assumptions of Definition \ref{definition.I.w.case.m.2}, assume also that all of the following
conditions hold:
\begin{enumerate}
\item[C1:] $\frac{\sigma}{\gamma} \le \frac{y}{q-1}$,
\item[C2.1:]
\be
 (2x_1+x_2) \left( y+\frac{\sigma}{\gamma}+2x_1+4x_2  \right)^2
 <  \left( 1-y -\frac{\sigma}{\gamma}-2(2x_1+3x_2) \right)^2 \left( y + \frac{\sigma}{\gamma}\right),
\nn\ee
\item[C2.2:]
\be
 2x_2 \left( y+\frac{\sigma}{\gamma}+2x_1+4x_2  \right)^4
< \left( 1-y -\frac{\sigma}{\gamma}- 2(2x_1+3x_2) \right)^3 \left( y + \frac{\sigma}{\gamma}\right)^2,
\nn\ee

\item[C3:] $\frac{\sigma}{\gamma} (1-y) < y^2$,

\item[C4:]
\be
x_2^2\left( y+ \frac{\sigma}{\gamma} + 2x_1 + 4 x_2 \right) \le 2x_1^3.
\nn\ee
\end{enumerate}
Then we have
$I_{y,x_1,x_2}(\sigma)=S(\sigma,y,0,2x_1,2x_2)$.
\end{proposition}
\begin{proof}
As in the proof of Corollary \ref{corollary.I.in.restricted.case}, we first observe that
for $0 \le x \le \frac{\sigma}{\gamma}$ and $0 \le t_1, t_2$ with $t_2 \le 2x_2$, $t_1 \le 2x_1 + x_2$, and
$2t_1 + 3t_2 \le 2(2x_1+3x_2)$, using condition C1 we obtain
\be \label{ep1.proposition.I.in.restricted.case.case.m.2}
\frac{y+x+t_1 + 2t_2}{y+\frac{\sigma}{\gamma} + t_1 + 2t_2} \ge \frac{y}{y+ \frac{\sigma}{\gamma}} \ge 1- \frac{1}{q}.
\ee

For the partial derivative $\frac{\partial S}{\partial x}(\sigma,y,x,t_1,t_2)$ of $S(\sigma,y,x,t_1,t_2)$ with respect to $x$,
using  (\ref{ep1.proposition.I.in.restricted.case.case.m.2}) and some straightforward manipulations
we get
\be
\frac{\partial S}{\partial x}(\sigma,y,x,t_1,t_2)=\log_q \frac{ (1-y-x-2t_1-3t_2)\left(\frac{\sigma}{\gamma}-x\right)}{(y+x+t_1+2t_2)^2}.
\nn\ee
By condition C3 we have
\be
\left(\frac{\sigma}{\gamma} -x \right) (1-y-x-2t_1-3t_2) \le \frac{\sigma}{\gamma} (1-y)
< y^2 \le (y+x+t_1+2t_2)^2,
\nn\ee
and hence
\be
\frac{\partial S}{\partial x}(\sigma,y,x,t_1,t_2) \le \log_q \frac{\frac{\sigma}{\gamma}(1-y)}{y^2} < 0
\nn\ee
for $0 < x < \frac{\sigma}{\gamma}$ and $0 \le t_1, t_2$ with $t_2 \le 2x_2$, $t_1 \le 2x_1+x_2$, and $2t_1+3t_2 \le 2(2x_1+3x_2)$.

Now we assume that
\be \label{ep1.in.the.middle.2.proposition.I.in.restricted.case.case.m.2}
x_1 > 0 \ \ \mbox{and} \ \  x_2>0.
\ee
For the partial derivatives $\frac{\partial S}{\partial t_1}(\sigma,y,x,t_1,t_2)$ and
$\frac{\partial S}{\partial t_2}(\sigma,y,x,t_1,t_2)$ of $S(\sigma,y,x,t_1,t_2)$ with respect to $t_1$ and $t_2$, again
using  (\ref{ep1.proposition.I.in.restricted.case.case.m.2}) and some straightforward manipulations
we get
\be
\frac{\partial S}{\partial t_1}=\log_q \frac{ (1-y-x-2t_1-3t_2)^2\left(y+\frac{\sigma}{\gamma}+t_1+2t_2\right)}{(y+x+t_1+2t_2)^2t_1}
\nn\ee
and
\be
\frac{\partial S}{\partial t_2}=\log_q \frac{ (1-y-x-2t_1-3t_2)^3\left(y+\frac{\sigma}{\gamma}+t_1+2t_2\right)^2}{(y+x+t_1+2t_2)^4t_2}.
\nn\ee

Note that $t_1+2t_2$ assumes its maximum over the region
\be \label{ep2.proposition.I.in.restricted.case.case.m.2}
0 \le t_2 \le 2x_2, \; 0 \le t_1 \le 2x_1+x_2, \; \mbox{and} \; 2t_1+3t_2 \le 2(2x_1+3x_2)
\ee
when $t_1=2x_1$ and $t_2=2x_2$. Therefore we have
\be \label{ep3.proposition.I.in.restricted.case.case.m.2}
t_1+2t_2 \le 2x_1 + 4x_2
\ee
over the region (\ref{ep2.proposition.I.in.restricted.case.case.m.2}).

Using (\ref{ep3.proposition.I.in.restricted.case.case.m.2}) and condition C2.1, we obtain
\be
\begin{array}{rl}
t_1(y+x+t_1+2t_2)^2 & \le

(2x_1+x_2) \left( y+ \frac{\sigma}{\gamma} + 2x_1+4x_2 \right)^2 \\ \\
& <

\left( y+ \frac{\sigma}{\gamma} \right) \left( 1-y-\frac{\sigma}{\gamma}-2(2x_1+3x_2) \right)^2 \\ \\

& \le \left( y+\frac{\sigma}{\gamma}+t_1 + 2t_2 \right) \left( 1-y-x-2t_1-3t_2 \right)^2.
\end{array}
\nn\ee
Similarly, using (\ref{ep3.proposition.I.in.restricted.case.case.m.2}) and condition C2.2 we obtain
\be
\begin{array}{rl}
t_2(y+x+t_1+2t_2)^4 & \le

2x_2 \left( y + \frac{\sigma}{\gamma}  + 2x_1 + 4x_2 \right)^4   \\ \\

& < \left( y+ \frac{\sigma}{\gamma} \right)^2 \left( 1-y-\frac{\sigma}{\gamma}-2(2x_1+3x_2) \right)^3 \\ \\

& \le \left(y+\frac{\sigma}{\gamma}+t_1 + 2t_2\right)^2 \left( 1-y-x-2t_1-3t_2 \right)^3.
\end{array}
\nn\ee
Hence we have
\be
\frac{\partial S}{\partial t_1} \ge \frac{ \left(1-y-\frac{\sigma}{\gamma}-2(2x_1+3x_2)\right)^2 \left(y+\frac{\sigma}{\gamma}\right)}{(2x_1+x_2) \left(y+\frac{\sigma}{\gamma}+2x_1+4x_2\right)^2} > 0
\nn\ee
and
\be
\frac{\partial S}{\partial t_2} \ge \frac{ \left(1-y-\frac{\sigma}{\gamma}-2(2x_1+3x_2)\right)^3 \left(y+\frac{\sigma}{\gamma}\right)^2}{2x_2 \left(y+ \frac{\sigma}{\gamma} + 2x_1+4x_2\right)^4} > 0
\nn\ee
for $0 \le x \le \frac{\sigma}{\gamma}$ and $0< t_1,t_2$ with $ t_2 \le 2x_2$, $t_1 \le 2x_1+x_2$, and $2t_1+3t_2 \le 2(2x_1+3x_2)$.

Then for fixed $\sigma$, $y$, and $0 \le x \le \frac{\sigma}{\gamma}$, the
function $S(\sigma,y,x,t_1,t_2)$ assumes its maximum
over the region (\ref{ep2.proposition.I.in.restricted.case.case.m.2}) on the part of the boundary formed by the closed line
connecting the two points
\be
A_1=(2x_1,2x_2) \ \ \mbox{and} \ \ A_2=\left(2x_1+x_2, \frac{4}{3}x_2\right).
\nn\ee
The direction vector $\overrightarrow{A_2A_1}$ from $A_2$ to $A_1$ is parallel to the vector $(-3,2)$. Hence
for  fixed $\sigma$, $y$, and $0 \le x \le \frac{\sigma}{\gamma}$,
the function $S(\sigma,y,x,t_1,t_2)$ is nondecreasing on the
closed line
from $A_2$ to $A_1$ if
\be \label{ep6.proposition.I.in.restricted.case.case.m.2}
-3 \frac{\partial S}{\partial t_1} (\sigma,y,x,t_1,t_2)
+2 \frac{\partial S}{\partial t_2} (\sigma,y,x,t_1,t_2) \ge 0
\ee
holds for fixed $\sigma$, $y$, and $0 \le x \le \frac{\sigma}{\gamma}$ and for each point $(t_1,t_2)$ on the
closed line from $A_2$ to $A_1$.
By straightforward manipulations, we obtain that (\ref{ep6.proposition.I.in.restricted.case.case.m.2})
is equivalent to
\be \label{ep7.proposition.I.in.restricted.case.case.m.2}
t_1^3\left( y + \frac{\sigma}{\gamma} + t_1+ 2t_2\right) \ge t_2^2 (y+x+t_1+2t_2)^2.
\ee
We have $t_1 \ge 2x_1$, $t_2 \le 2x_2$,  and $t_1+2t_2 \le 2x_1+4x_2$
on the closed line from $A_2$ to $A_1$.
Therefore using
$ y + \frac{\sigma}{\gamma} + t_1+ 2t_2  \ge y+x+t_1+2t_2$ and condition C4,
we see that (\ref{ep7.proposition.I.in.restricted.case.case.m.2})
holds. Hence $S(y,\sigma,x,t_1,t_2)$ assumes its maximum at $x=0$ and $(t_1,t_2)=A_1=(2x_1,2x_2)$. It is easy to check that
if the assumption (\ref{ep1.in.the.middle.2.proposition.I.in.restricted.case.case.m.2}) does not hold, but the assumptions of the proposition hold,
then similar methods also apply
and we again have $I_{t,x_1,x_2}(\sigma)=S(\sigma,y,0,2x_1,2x_2)$.
This completes the proof.
\end{proof}

Now that we have dealt with the cases $m=1$ and $m=2$, we present the generalizations for any $m \ge 1$.

\begin{definition} \label{definition.S.general.m}
Under Assumption 1 (cf. Section \ref{section.asymptotic.size.V.m.case.m.1}),
let $y>0$, $x_1,x_2, \ldots, x_m,\sigma \ge 0$ be real numbers satisfying
\be \label{condition.general.y.x1.x2.xm}
y+2(2x_1+3x_2+ \cdots+(m+1)x_m) + \frac{\sigma}{\gamma} < 1.
\ee
For real numbers $0 \le x \le \frac{\sigma}{\gamma}$ and $t_1, t_2, \ldots, t_m$ satisfying
\be \label{condition1.general.region.t1.t2.tm}
\begin{array}{l}
0 \le t_m \le 2x_m, \; 0 \le t_{m-1} \le 2x_{m-1} + x_m, \ldots , \\
0 \le t_1 \le 2x_1 + (x_2+x_3+ \cdots + x_m),
\end{array}
\ee
and
\be \label{condition2.general.region.t1.t2.tm}
2t_1+3t_2+ \cdots + (m+1)t_m \le 2(2x_1+3x_2+ \cdots+ (m+1)x_m),
\ee
let $S(\sigma,y,x,t_1,t_2, \ldots, t_m)$ be the real-valued function
\be
\begin{array}{l}
\dd S(\sigma,y,x,t_1,t_2, \ldots, t_m) \\ \\
=  E(t_m) + (1-t_m)E\left( \frac{t_{m-1}}{1-t_m} \right) + \cdots +
\left(1-(t_2+ \cdots+ t_m)\right) E \left( \frac{t_1}{1-(t_2+\cdots+t_m)} \right) \\ \\

\dd
+ \left(1-(t_1+t_2+ \cdots + t_m)\right)E\left( \frac{y+x+(t_1+2t_2+ \cdots + mt_m)}{1-(t_1+t_2+ \cdots + t_m)}\right) \\ \\

+ \left\{
\begin{array}{l}
\left(y+ \frac{\sigma}{\gamma} + (t_1+2t_2+ \cdots + mt_m) \right)
E\left( \frac{y+x+(t_1+2t_2+ \cdots + mt_m)}{y+ \frac{\sigma}{\gamma} + (t_1+2t_2+ \cdots + mt_m)}\right) \\ \\
\hspace{5cm} \mbox{if $ \frac{y+x+(t_1+2t_2+ \cdots + mt_m)}{y+\frac{\sigma}{\gamma} + (t_1+2t_2+ \cdots + mt_m)} \ge 1 -\frac{1}{q}$}, \\ \\

\left(y+\frac{\sigma}{\gamma}+(t_1+2t_2+ \cdots +mt_m)\right) - (y+x+(t_1+2t_2+ \cdots +mt_m))\log_q (q-1) \\ \\
\hspace{5cm} \mbox{if $ \frac{y+x+(t_1+2t_2+ \cdots + mt_m)}{y+\frac{\sigma}{\gamma} + (t_1+2t_2+ \cdots + mt_m)} \le 1 -\frac{1}{q}$}.
\end{array}
\right.
\end{array}
\nn\ee
Note that by (\ref{condition.general.y.x1.x2.xm}) we have $2\left(2x_1+3x_2+ \cdots (m+1)x_m\right) <1$, and hence using
(\ref{condition2.general.region.t1.t2.tm}) we obtain
$t_1 +t_2 + \cdots +t_m \le t_1 + \frac{3}{2} t_2 + \cdots + \frac{m+1}{2}t_m < \frac{1}{2}$.
\end{definition}

We state the generalization of Proposition \ref{propostion.asymptotic.card.V.m.case.2} whose proof is similar.

\begin{proposition} \label{general.case.propostion.asymptotic.card.V.m}
Under Assumption 1 (cf. Section \ref{section.asymptotic.size.V.m.case.m.1}), let
$y > 0$ and $x_1,x_2,\ldots,x_m,\sigma \ge 0$ be real numbers satisfying (\ref{condition.general.y.x1.x2.xm}).
For each integer $i \ge 1$,
let $r_i=\left\lfloor \left(m + y + \frac{\sigma}{\gamma} \right) n_i \right \rfloor$,
$s_i=\left\lfloor (m+y)n_i \right \rfloor$, $X_1^{(i)}=\left \lfloor x_1n_i \right \rfloor$,
$X_2^{(i)}=\left \lfloor x_2 n_i \right\rfloor$, \ldots,
$X_m^{(i)}=\left \lfloor x_m n_i \right\rfloor$,
and
$\mathcal{V}^{(i)}_m(r_i,s_i;X_1^{(i)},X_2^{(i)},\ldots, X_m^{(i)})$ be the set of positive divisors of $F_i$
defined in Definition \ref{definition.V.m}. Then for the cardinalities of these sets we have
\be
\limsup_{i \ra \infty} \frac{ \log_q \left| \mathcal{V}_m^{(i)}(r_i,s_i;X_1^{(i)},X_2^{(i)}, \ldots, X_m^{(i)})\right|}{n_i}
\le \max S(\sigma,y,x,t_1,t_2, \ldots, t_m),
\nn\ee
where the maximum is over all real numbers $x$ and $t_1, t_2, \ldots, t_m$ satisfying $0 \le x \le \frac{\sigma}{\gamma}$
and  the conditions in (\ref{condition1.general.region.t1.t2.tm}) and (\ref{condition2.general.region.t1.t2.tm}).
\end{proposition}

Now we generalize Definition \ref{definition.I.w.case.m.2}.

\begin{definition} \label{general.definition.I.y.x1.x2.xm}
Under Assumption 1 (cf. Section \ref{section.asymptotic.size.V.m.case.m.1}), let $y>0$ and $x_1,x_2, \ldots, x_m\ge 0$
be real numbers
such that $y+2\left(2x_1+3x_2+ \cdots + (m+1)x_m\right) < 1$. For $\sigma \ge 0$ and
$y+2\left(2x_1+3x_2+ \cdots + (m+1)x_m\right)+\frac{\sigma}{\gamma} <1$, let
$I_{y,x_1,x_2, \ldots, x_m}(\sigma)$ be the real-valued function of $\sigma$ defined by
\be
I_{y,x_1,x_2, \ldots, x_m}(\sigma)=\max S(\sigma,y,x,t_1,t_2, \ldots, t_m),
\nn\ee
where the maximum is over all real numbers $x,t_1,t_2, \ldots, t_m$ with $0 \le x \le \frac{\sigma}{\gamma}$ and
$t_1, t_2, \ldots, t_m$ satisfying conditions (\ref{condition1.general.region.t1.t2.tm}) and (\ref{condition2.general.region.t1.t2.tm}).
\end{definition}

The proof of the next lemma generalizing Lemma \ref{lemma.I.increasing.case.m.2} is also similar.

\begin{lemma} \label{general.lemma.I.increasing}
Under the assumptions of Definition \ref{general.definition.I.y.x1.x2.xm},
the real-valued function $I_{y,x_1,x_2, \ldots, x_m}(\sigma)$ is a strictly
increasing function of $\sigma$ on its domain of definition, which is the interval of $\sigma$ such that
$\sigma \ge 0$ and $y+2\left(2x_1+3x_2+ \cdots +(m+1)x_m\right)+\frac{\sigma}{\gamma} < 1$.
\end{lemma}

Now we are ready to compute $I_{y,x_1,x_2, \ldots, x_m}(\sigma)$ for general $m$ under some conditions.
We note that since the region defined by the conditions (\ref{condition1.general.region.t1.t2.tm}) and (\ref{condition2.general.region.t1.t2.tm})
is more complicated in the general case than the one in the case $m=2$, we need to define new parameters in
the following proposition in order to state the result.

\begin{proposition} \label{general.proposition.I.y.x1.x2.xm.compute}
Under the assumptions of Definition \ref{general.definition.I.y.x1.x2.xm}, let
\be
\bar{t}_m=2x_m \ \ \mbox{and} \ \  \bar{t}_\ell=2x_\ell + \sum_{\nu = \ell +1}^m x_\nu \hspace{1cm} \mbox{for} \; 1 \le \ell \le m-1.
\nn\ee
Let $t_1^*$ be the real number defined by
\be
2t_1^*+\sum_{\ell =2}^m (\ell +1) \bar{t}_\ell = 2 \sum_{\ell=1}^m (\ell +1) x_\ell,
\nn\ee
and for each $2 \le \ell \le m$, let $t_\ell^*$ be the real number defined inductively using $t_{\ell-1}^*$ by
\be \label{e1.general.proposition.I.y.x1.x2.xm.compute}
(\ell+1) t_\ell^* - (\ell +1) \bar{t}_\ell = \ell t_{\ell-1}^* - \ell \bar{t}_{\ell -1}.
\ee
Moreover, let $u$ be the real number depending on $x_1, \ldots, x_m$ defined by
\be
u=t_1^* + \sum_{\ell=2}^m \ell \bar{t}_\ell.
\nn\ee
Assume also that all of the following
conditions hold:
\begin{enumerate}
\item[C1:] $\frac{\sigma}{\gamma} \le \frac{y}{q-1}$,
\item[C2:] For each $1 \le \ell \le m$,
\be
\bar{t}_\ell \left( y + \frac{\sigma}{\gamma} + u\right)^{2\ell} <
\left(1 - y - \frac{\sigma}{\gamma} -2\sum_{\nu=1}^m (\nu+1)x_\nu\right)^{\ell+1} \left( y + \frac{\sigma}{\gamma} \right)^\ell,
\nn\ee

\item[C3:] $\frac{\sigma}{\gamma} (1-y) < y^2$,

\item[C4:] For each $1 \le \ell \le m-1$,
\be
\left(\bar{t}_{\ell+1}\right)^{\ell+1} \left( y + \frac{\sigma}{\gamma} + u \right) \le \left(t_\ell^*\right)^{\ell+2}.
\nn\ee
\end{enumerate}
Then we have $I_{y,x_1,x_2,x_3, \ldots, x_m}(\sigma)=S(\sigma,y,0,t_1^*,\bar{t}_2,\bar{t}_3, \ldots, \bar{t}_m)$.
\end{proposition}
\begin{proof}
By condition C1 we have
\be \label{ep1.general.proposition.I.y.x1.x2.xm.compute}
\frac{y+x+t_1 + 2t_2+ \cdots +mt_m}{y+\frac{\sigma}{\gamma} + t_1 + 2t_2+ \cdots + mt_m}
\ge \frac{y}{y+ \frac{\sigma}{\gamma}} \ge 1- \frac{1}{q}.
\ee
The following identities for partial derivatives hold:
\be
\begin{array}{l}
\dd \frac{\partial}{\partial x} \left\{
\left(1-t_1-t_2- \cdots - t_m \right)
E \left( \frac{y+x+t_1+2t_2+ \cdots + mt_m}{1-t_1-t_2-\cdots-t_m}\right) \right\}\\
\dd =\log_q \frac{1-y-x-2t_1-3t_2 - \cdots -(m+1)t_m}{y+x+t_1+2t_2+ \cdots + mt_m}
\end{array}
\nn\ee
and
\be
\begin{array}{l}
\dd \frac{\partial}{\partial x}\left\{
\left( y + \frac{\sigma}{\gamma} + t_1+ 2t_2 + \cdots + mt_m \right)
E\left( \frac{y+x+t_1+2t_2+ \cdots + mt_m}{y+\frac{\sigma}{\gamma} + t_1+ 2t_2 + \cdots + mt_m}\right)\right\} \\
\dd = \log_q \frac{\frac{\sigma}{\gamma}-x}{y+x+t_1+ 2t_2+ \cdots + mt_m}.
\end{array}
\nn\ee
Hence using Definition \ref{definition.S.general.m} and  (\ref{ep1.general.proposition.I.y.x1.x2.xm.compute}), we obtain that
\be
\frac{\partial S}{\partial x}
= \log_q \frac{(1-y-x-2t_1-3t_2-\cdots-(m+1)t_m)\left(\frac{\sigma}{\gamma}-x\right)}{(y+x+t_1+ 2t_2 + \cdots+mt_m)^2}.
\nn\ee
Therefore if conditions C1 and C3 hold, then
\be \label{ep2.general.proposition.I.y.x1.x2.xm.compute}
\frac{\partial S}{\partial x}(\sigma,y,x,t_1, t_2,\ldots,t_m) <0
\ee
for each
$0 < x < \frac{\sigma}{\gamma}$ and $t_1, \ldots, t_m$ in the region defined by
(\ref{condition1.general.region.t1.t2.tm}) and (\ref{condition2.general.region.t1.t2.tm}).

Now we further assume that
\be \label{ep3.general.proposition.I.y.x1.x2.xm.compute}
x_1 > 0, \; x_2> 0, \ldots, x_m > 0.
\ee
For $1 \le \ell \le m$, by straightforward manipulations we also obtain the following identities for  partial derivatives:
\be
\begin{array}{l}
 \frac{\partial}{\partial t_\ell} \left\{
E(t_m) + (1-t_m) E \left( \frac{t_{m-1}}{1-t_m}\right) + \cdots + (1-t_2- \cdots-t_m) E \left( \frac{t_1}{1-t_2-\cdots-t_m}\right) \right\} \\
\dd = \log_q(1-t_1-t_2 - \cdots -t_m) - \log_q(t_\ell),
\end{array}
\nn\ee
\be
\begin{array}{rl}
& \frac{\partial}{\partial t_\ell} \left\{
(1-t_1-t_2 -\cdots -t_m) E \left( \frac{y+x+ t_1+2t_2+ \cdots + mt_m}{1-t_1-t_2 -\cdots -t_m} \right) \right\} \\
=& \dd(\ell +1) \log_q \left( 1-y-x-2t_1-3t_2 -\cdots -(m+1)t_m\right) \\
&\dd -\log_q(1-t_1-t_2-\cdots-t_m) -\ell \log_q(y+x+t_1+2t_2+ \cdots +mt_m),
\end{array}
\nn\ee
and
\be
\begin{array}{l}
\frac{\partial}{\partial t_\ell} \left\{
\left(y+\frac{\sigma}{\gamma} + t_1+ 2t_2 + \cdots + mt_m \right)
E \left( \frac{y+x+t_1+2t_2+ \cdots + mt_m}{y+ \frac{\sigma}{\gamma} +t_1+ 2t_2+ \cdots +mt_m}\right)\right\} \\
\dd =
\ell \log_q\left(y+ \frac{\sigma}{\gamma} + t_1+2t_2 + \cdots + mt_m\right)
-\ell \log_q (y+x+t_1+2t_2+ \cdots + mt_m).
\end{array}
\nn\ee
Hence using Definition \ref{definition.S.general.m} and (\ref{ep1.general.proposition.I.y.x1.x2.xm.compute}), for $1 \le \ell \le m$ we obtain that
\be \label{ep4.general.proposition.I.y.x1.x2.xm.compute}
\begin{array}{rl}
 \dd \frac{\partial S}{\partial t_\ell}=
& \dd \log_q (1-y-x-2t_1-3t_2 -\cdots -(m+1)t_m)^{\ell +1} \\
& \dd + \log_q \left(y+ \frac{\sigma}{\gamma} +t_1+ 2t_2+ \cdots + mt_m\right)^\ell \\
& - \log_q (y+x+t_1+2t_2+ \cdots +mt_m)^{2\ell} - \log_q t_\ell.
\end{array}
\ee
Now we also assume that for the real number $u$ defined in the statement of the proposition we have
\be \label{ep5.general.proposition.I.y.x1.x2.xm.compute}
u=\max \left( t_1 + 2t_2 + \cdots + mt_m\right),
\ee
where the maximum is over the region defined by the conditions (\ref{condition1.general.region.t1.t2.tm}) and (\ref{condition2.general.region.t1.t2.tm}).
Later in this proof, we will show that the assumption (\ref{ep5.general.proposition.I.y.x1.x2.xm.compute}) holds.

Using (\ref{ep1.general.proposition.I.y.x1.x2.xm.compute}), (\ref{ep4.general.proposition.I.y.x1.x2.xm.compute}),
(\ref{ep5.general.proposition.I.y.x1.x2.xm.compute}), and condition C2,
as in the proof of Proposition \ref{proposition.I.in.restricted.case.case.m.2}, we obtain that for each $1 \le \ell \le m$,
\be
\frac{\partial S}{\partial t_\ell} \left(\sigma,y,x,t_1, \ldots, t_m\right) >0
\nn\ee
holds for  $ 0 \le x \le \frac{\sigma}{\gamma}$ and the real numbers $0< t_1, \ldots, t_m$ satisfying the conditions
(\ref{condition1.general.region.t1.t2.tm}) and (\ref{condition2.general.region.t1.t2.tm}). This implies that
for each $0 \le x \le \frac{\sigma}{\gamma}$, $S(\sigma,y,x,t_1, \ldots, t_m)$ assumes its maximum over the region
defined by (\ref{condition1.general.region.t1.t2.tm}) and (\ref{condition2.general.region.t1.t2.tm}) on the closed set, forming a part of the
boundary of the region, defined by the conditions
\be \label{ep6.general.proposition.I.y.x1.x2.xm.compute}
0 \le t_\ell \le \bar{t}_\ell \hspace{1cm} \mbox{for} \; 1 \le \ell \le m
\ee
and
\be \label{ep7.general.proposition.I.y.x1.x2.xm.compute}
\sum_{\ell=1}^m (\ell+1) t_\ell = 2 \sum_{\ell=1}^m (\ell +1) x_\ell,
\ee
where $\bar{t}_\ell$ is defined in the statement of the proposition.

For each $1 \le \ell \le m$, it follows from the definition of $t_\ell^*$ in the statement of the proposition that
$t_\ell^*$ is the smallest value of the parameter $t_\ell$ over the closed set defined by the conditions (\ref{ep6.general.proposition.I.y.x1.x2.xm.compute})
and (\ref{ep7.general.proposition.I.y.x1.x2.xm.compute}). For each $1 \le \ell \le m$, let $A_\ell$ be the point of the
$(t_1, \ldots, t_m)$-space given by
\be
A_\ell=(t_1, \ldots, t_m) \hspace{0.7cm} \mbox{where} \; t_\ell=t_\ell^* \; \mbox{and} \; t_\nu=\bar{t}_\nu \; \mbox{for} \; \nu \in \{1, \ldots, m\} \setminus \{\ell\}.
\nn\ee
We observe that the points $A_1,A_2, \ldots, A_m$ are the corners of the closed set given by (\ref{ep6.general.proposition.I.y.x1.x2.xm.compute}) and (\ref{ep7.general.proposition.I.y.x1.x2.xm.compute}).

For each $1 \le \ell \le m-1$, the direction vector $\overrightarrow{A_{\ell+1}A_\ell}$ from $A_{\ell+1}$ to $A_\ell$
in the $(t_1, \ldots, t_m)$-space is
\be
\overrightarrow{A_{\ell+1}A_\ell}
=\left(\underbrace{0,\ldots,0}_{\ell -1 \;
\mbox{times}},t_\ell^*-\bar{t}_\ell, \bar{t}_{\ell+1}-t_{\ell+1}^*, \underbrace{0, \ldots, 0}_{m-\ell-1\;
\mbox{times}}\right).
\nn\ee
Using (\ref{e1.general.proposition.I.y.x1.x2.xm.compute}) we observe that for each $1 \le \ell \le m-1$, the direction vector
$\overrightarrow{A_{\ell+1}A_\ell}$ is parallel to the vector
\be \label{ep8.general.proposition.I.y.x1.x2.xm.compute}
\left(\underbrace{0, \ldots, 0}_{\ell -1 \; \mbox{times}}, -(\ell+2),\ell+1,\underbrace{0,\ldots,0}_{m-\ell-1\;\mbox{times}} \right)
\ee
in the $(t_1, \ldots, t_m)$-space.

If for each $1 \le \ell \le m-1$ the inequality
\be \label{ep9.general.proposition.I.y.x1.x2.xm.compute}
\overrightarrow{A_{\ell+1}A_\ell}
\cdot
\left( \frac{\partial S}{\partial t_1}, \ldots,
\frac{\partial S}{\partial t_m}\right) (\sigma, y,x, t_1, \ldots, t_m) \ge 0
\ee
for the standard inner product of vectors in the $(t_1, \ldots,t_m)$-space holds for each $0 \le x \le \frac{\sigma}{\gamma}$
and $t_1, \ldots, t_m$ satisfying (\ref{ep6.general.proposition.I.y.x1.x2.xm.compute}) and (\ref{ep7.general.proposition.I.y.x1.x2.xm.compute}),
then $S(\sigma,y,x,t_1, \ldots, t_m)$ is nondecreasing in the directions from $A_m$ to $A_{m-1}$, from $A_{m-1}$ to $A_{m-2}$, \ldots, and
from $A_2$ to $A_1$. This implies that if (\ref{ep9.general.proposition.I.y.x1.x2.xm.compute}) holds, then for each
$0 \le x \le \frac{\sigma}{\gamma}$, $S(\sigma,y,x,t_1, \ldots, t_m)$ assumes its maximum at $A_1$. Using (\ref{ep8.general.proposition.I.y.x1.x2.xm.compute}),
we obtain that (\ref{ep9.general.proposition.I.y.x1.x2.xm.compute}) is equivalent to
\be \label{ep10.general.proposition.I.y.x1.x2.xm.compute}
(\ell+1) \frac{\partial S}{\partial t_{\ell+1}} (\sigma,y,x, t_1, \ldots, t_m) \ge (\ell +2) \frac{\partial S}{\partial t_\ell} (\sigma,y,x,t_1, \ldots, t_m).
\ee
Using (\ref{ep4.general.proposition.I.y.x1.x2.xm.compute}), (\ref{ep10.general.proposition.I.y.x1.x2.xm.compute}),
and some straightforward
manipulations, we observe that (\ref{ep9.general.proposition.I.y.x1.x2.xm.compute}) holds if
\be \label{ep11.general.proposition.I.y.x1.x2.xm.compute}
\begin{array}{l}
\left(t_\ell\right)^{\ell+2} \left(y + \frac{\sigma}{\gamma} + t_1+2t_2+ \cdots + mt_m\right) \\ \\
\ge
\left( t_{\ell+1}\right)^{\ell+1} (y+x+t_1+2t_2+ \cdots+ mt_m)^2.
\end{array}
\ee
Using the fact that $y + \frac{\sigma}{\gamma} + t_1+ 2t_2 + \cdots + mt_m \ge y+x+ t_1+ 2t_2+ \cdots +mt_m$,
the assumption (\ref{ep5.general.proposition.I.y.x1.x2.xm.compute}),
and condition C4, as in the proof of Proposition \ref{proposition.I.in.restricted.case.case.m.2}, we obtain that (\ref{ep11.general.proposition.I.y.x1.x2.xm.compute})
holds, and hence for each $0 \le x \le \frac{\sigma}{\gamma}$, $S(\sigma,y,x,t_1, \ldots, t_m)$ assumes its maximum at $A_1$.

Next we prove the claim (\ref{ep5.general.proposition.I.y.x1.x2.xm.compute}). Note that the gradient of the $m$-variable function
$f(t_1,t_2, \ldots, t_m)=t_1+2t_2+ \cdots + mt_m$ is $(1,2, \ldots, m)$ at any point of the $(t_1, \ldots, t_m)$-space.
For each $1 \le \ell \le m-1$, from the standard inner product with the vector in (\ref{ep8.general.proposition.I.y.x1.x2.xm.compute})
we obtain
\be
\begin{array}{l}
(1,2, \ldots,m) \cdot \left(\underbrace{0, \ldots, 0}_{\ell -1 \; \mbox{times}}, -(\ell+2),\ell+1,\underbrace{0,\ldots,0}_{m-\ell-1\;\mbox{times}} \right) \\ \\
= -(\ell +2) \ell + (\ell+1)^2 =1 >0.
\end{array}
\nn\ee
Then, as the function $S(\sigma,y,x,t_1, \ldots, t_m)$, the function $f(t_1, \ldots, t_m)$ assumes its maximum at $A_1$
and hence the claim (\ref{ep5.general.proposition.I.y.x1.x2.xm.compute}) holds. Finally, using (\ref{ep2.general.proposition.I.y.x1.x2.xm.compute})
we complete the proof under the assumption (\ref{ep3.general.proposition.I.y.x1.x2.xm.compute}).
As in the proof of Proposition \ref{proposition.I.in.restricted.case.case.m.2}, we observe that
if the assumption (\ref{ep3.general.proposition.I.y.x1.x2.xm.compute}) does not hold, but the assumptions of the proposition hold,
then similar methods also apply
and we again have $I_{y,x_1,x_2, x_3,\ldots, x_m}(\sigma)=S(\sigma,y,0,t_1^*,\bar{t}_2, \bar{t}_3, \ldots, \bar{t}_m)$.
This completes the proof.
\end{proof}
\begin{remark} \label{general.remark.proposition.I.y.x1.x2.xm.compute}
We note that Proposition \ref{general.proposition.I.y.x1.x2.xm.compute} reduces to Proposition \ref{proposition.I.in.restricted.case.case.m.2}
and Corollary \ref{corollary.I.in.restricted.case} if $m=2$ and $m=1$, respectively.
\end{remark}

\section{Asymptotic Bounds for Codes} \label{section.asymptotic.constructon}

In this section we prove our main results (Theorem \ref{theorem.asymp} and Corollary \ref{corollary.asym.linear})
which establish improved lower bounds on $\alpha_q(\delta)$ and $\alpha_q^{\rm lin}(\delta)$.

We first state our main assumption, which is like Assumption 1 in Section \ref{section.asymptotic.size.V.m.case.m.1},
but introduces more notation.

\begin{description}

 \item[{\bf Assumption $\boldsymbol{1^\prime}$}] Assume that $\left(F_i/\F_q\right)_{i=1}^\infty$ is a sequence of global function fields with full
constant field $\F_q$, with $g_i \ra \infty$
as $i \ra \infty$, and with $\lim_{i \ra \infty} \frac{n_i}{g_i}=\gamma >0$, where $n_i$ and $g_i$ denote
the number of rational places and the genus of $F_i$,
respectively. For each $l \ge 1$, let $\gamma_l \ge 0$ be a real number
with $\liminf_{i \ra \infty} \dd \frac{B_{i,l}}{g_i} \ge \gamma_l$,
where $B_{i,l}$ is the number of degree $l$ places of $F_i$. Note that
we can take $\gamma_1=\gamma$.

\end{description}

The following well-known result will be useful.

\begin{proposition} \label{proposition.class.number.upper.bound}
Under Assumption $1^\prime$
we have
\be
\liminf_{i \ra \infty} \frac{\log_q h_i}{n_i} \ge \frac{1}{\gamma} \left[ 1 + \sum_{l=1}^\infty \gamma_l \log_q \frac{q^l}{q^l-1} \right],
\nn\ee
where $h_i$ is the class number of $F_i$.
\end{proposition}
\begin{proof}
This follows from \cite[Corollary 2]{T} (see also \cite[Exercise 2.3.27]{TV}).
\end{proof}

Now we introduce an important function based on the function $I_{y,x_1,\ldots,x_m}(\sigma)$ defined in Definition \ref{general.definition.I.y.x1.x2.xm}.
In the next definition we use the fact that $I_{y,x_1, \ldots, x_m}(\sigma)$ is an increasing function on its domain of definition, see Lemma \ref{general.lemma.I.increasing}.

\begin{definition} \label{definition.Psi}
Under Assumption $1^\prime$ and
for real numbers $y>0$ and $x_1, \ldots, x_m \ge 0$ with $y+2(2x_1+3x_2+ \cdots + (m+1)x_m) <1$,
let $\Psi(y,x_1, \ldots, x_m)$ be the real-valued function of $y,x_1, \ldots, x_m$ defined
by
\be
\Psi(y,x_1, \ldots, x_m)=
\left\{\begin{array}{l}
\dd I^{-1}_{y,x_1, \ldots, x_m}\left( \frac{1}{\gamma} \left[ 1 + \sum_{l=1}^\infty \gamma_l \log_q \frac{q^l}{q^l-1}\right] \right) \\ \\
\hspace{1.58cm} \mbox{if} \; \dd \lim_{\sigma \ra \theta^-}
I_{y,x_1, \ldots, x_m}(\sigma) >  \frac{1}{\gamma} \left[ 1 + \sum_{l=1}^\infty \gamma_l \log_q \frac{q^l}{q^l-1}\right],  \\
0
\hspace{1.4cm} \mbox{otherwise},
\end{array}
\right.
\nn\ee
where $\theta=\gamma\left(1-y-2(2x_1+3x_2 + \cdots (m+1)x_m)\right)$.
\end{definition}

Now we are ready to establish our main results. We recall that the functions $\alpha_q(\delta)$ and $\alpha_q^{\rm lin}(\delta)$
are defined in (\ref{definition.alpha.q.delta}) and (\ref{definition.alpha.linear.q.delta}), respectively.

\begin{theorem} \label{theorem.asymp}
Under Assumption $1^\prime$,
let $x_1, \ldots, x_m \ge 0$ be real numbers with
$2(2x_1+3x_2+ \cdots + (m+1)x_m) < 1$. For each real number $0 < \delta < 1-2(2x_1+3x_2+ \cdots + (m+1)x_m)$
we have
\be
\begin{array}{l}
\dd \alpha_q(\delta) \ge  \dd R_{\{\gamma_l\},x_1, \ldots, x_m}(\delta):=1-\delta - \frac{1}{\gamma} + (x_1+ \cdots + x_m) \log_q (q-1) \\ \\
 \hspace{0.5cm} \dd -(x_1\log_q x_1+ \cdots + x_m\log_q x_m) - \left(1-(x_1+ \cdots + x_m)\right) \log_q \left(1-(x_1+ \cdots + x_m)\right) \\ \\
 \hspace{0.5cm}\dd -\left(4x_1+5x_2+ \cdots +(m+3)x_m\right) \\ \\
 \hspace{0.5cm} \dd + \frac{1}{\gamma} \Psi\Big(1-\delta-2(2x_1+3x_2+ \cdots (m+1)x_m),x_1,x_2,\ldots,x_m\Big).
\end{array}
\nn\ee
\end{theorem}
\begin{proof}
Let $y=1-\delta-2(2x_1+3x_2+\cdots + (m+1)x_m)$
and $\sigma=\Psi(y,x_1, \ldots,x_m)$. If $\sigma=0$, then
the theorem follows from \cite[Theorem 5.1]{NO1}. If $R_{\{\gamma_l\},x_1,\ldots, x_m}(\delta) \le 0$, then the
statement of the theorem is trivial. Therefore we can assume that $\sigma >0$ and $ R_{\{\gamma_l\},x_1,\ldots, x_m}(\delta) > 0$.
Let $0 <\epsilon < \sigma  $ be a real number such that
\be \label{ep-1.theorem.asymp}
\begin{array}{l}
\hspace{0.5cm} y+(x_1+ \cdots +x_m) \log_q (q-1) - (x_1\log_q x_1 + \cdots + x_m \log_q x_m) \\
\hspace{0.5cm} -\left(1-(x_1+ \cdots + x_m) \right) \log_q \left(1-(x_1+ \cdots + x_m) \right) \\
\hspace{0.5cm}+ (x_2+ 2x_3 + \cdots + (m-1)x_m)  \\ \\
> \dd \frac{1-(\sigma-\epsilon)}{\gamma}.
\end{array}
\ee
For $i \ge 1$, let
\be \label{ep0.theorem.asymp}
\begin{array}{l}
r_i=\left\lfloor \left(m+y+ \frac{\sigma -\epsilon}{\gamma}\right) n_i \right\rfloor,  s_i=\left\lfloor (m+y)n_i \right\rfloor, \\
X_1^{(i)}=\lfloor x_1 n_i \rfloor, X_2^{(i)}= \lfloor x_2 n_i \rfloor, \ldots, X_m^{(i)}= \lfloor x_m n_i \rfloor.
\end{array}
\ee
For sufficiently large $i$, by Propositions \ref{general.case.propostion.asymptotic.card.V.m} and
\ref{proposition.class.number.upper.bound},
the hypotheses of Proposition \ref{proposition.existence.G} for the global function field $F_i$
with $r_i$, $s_i$, and $X_1^{(i)}, \ldots, X_m^{(i)}$ as in (\ref{ep0.theorem.asymp}) are satisfied.
Let $G_i$ be the divisor of $F_i$ given by
Proposition \ref{proposition.existence.G} with these parameters for sufficiently large $i$.

Note that
\be \label{ep1.theorem.asymp}
\begin{array}{l}
\dd \liminf_{i \ra \infty} \frac{\log_q \left| M( x_1 , \ldots, x_m ; \boldsymbol{0})\right| }{n_i} \\
\ge (x_1+ \cdots + x_m) \log_q(q-1) - (x_1 \log_q x_1 + \cdots + x_m \log_q x_m) \\
\hspace{0.2cm} -\left( 1 - (x_1+\cdots+x_m)\right) \log_q \left( 1 - (x_1+\cdots+x_m)\right) \\
\hspace{0.2cm} + (x_2+2x_3+ \cdots+(m-1)x_m)
\end{array}
\ee
(see \cite[Section 4]{NO1}). Since we have (\ref{ep-1.theorem.asymp}),
using the divisor $G_i$ of the global
function field $F_i$ for sufficiently large $i$, Theorem \ref{theorem.basic}, and (\ref{ep1.theorem.asymp}),
we obtain a sequence of $q$-ary codes $\{C_i\}_{i=1}^\infty$ of lengths $\{n_i\}_{i=1}^\infty$, respectively, such that
$n_i \ra \infty$ as $i \ra \infty$ by Assumption $1^\prime$ as well as
\be
\begin{array}{rl}
&\dd \liminf_{i \ra \infty} \frac{\log_q |C_i|}{n_i}   \\
\ge & \dd  y + \frac{\sigma-\epsilon}{\gamma} - \frac{1}{\gamma}  \\
& \dd + (x_1+ \cdots +x_m) \log_q (q-1)  -(x_1\log_q x_1 + \cdots + x_m \log_q x_m) \\
& \dd - \left(1-(x_1+ \cdots + x_m)\right) \log_q\left(1-(x_1+ \cdots + x_m)\right) \\
& \dd +(x_2+2x_3+ \cdots + (m-1)x_m) \\
= & \dd 1-\delta-2(2x_1+3x_2+\cdots + (m+1)x_m) + \frac{\sigma-\epsilon}{\gamma} - \frac{1}{\gamma} \\
& \dd + (x_1+ \cdots +x_m) \log_q (q-1)  -(x_1\log_q x_1 + \cdots + x_m \log_q x_m) \\
& \dd - \left(1-(x_1+ \cdots + x_m)\right) \log_q\left(1-(x_1+ \cdots + x_m)\right) \\
& \dd +(x_2+2x_3+ \cdots + (m-1)x_m) \\
=& \dd R_{\{\gamma_l\},x_1, \ldots, x_m}(\delta) - \frac{\epsilon}{\gamma}
\end{array}
\nn\ee
and
\be
\dd \liminf_{i \ra \infty} \frac{d(C_i)}{n_i} \ge \delta.
\nn\ee
Using the fact that $\alpha_q(\delta)$ is a nonincreasing function
of $\delta$, we get
\be
\alpha_q(\delta) \ge R_{\{\gamma_l\},x_1, \ldots, x_m}(\delta) - \frac{\epsilon}{\gamma}.
\nn\ee
Letting $ \epsilon \ra 0^+$
completes the proof.
\end{proof}

\begin{corollary} \label{corollary.asym.linear}
Under Assumption $1^\prime$, for each real number $0 < \delta < 1$
we have
\be
\begin{array}{l}
\dd \alpha^{\rm lin}_q(\delta) \ge  \dd R^{\rm lin}_{\{\gamma_l\}}(\delta):=1-\delta - \frac{1}{\gamma}
  +\frac{1}{\gamma}\Psi\Big(1-\delta,0\Big).
\end{array}
\nn\ee
\end{corollary}
\begin{proof}
Taking $m=1$ and using similar methods as in the proof of Theorem \ref{theorem.asymp}, but applying
Corollary \ref{corollary.basic.linear} instead of Theorem \ref{theorem.basic}, we obtain the desired result.
\end{proof}

\section{Examples} \label{section.example}

In this section we demonstrate that Theorem \ref{theorem.asymp} and Corollary \ref{corollary.asym.linear} yield improvements
on the lower bounds for $\alpha_q(\delta)$ and $\alpha^{\rm lin}_q(\delta)$ at least for certain values of $q$ and certain
values of $\delta$. In our examples we use well-known values for $\gamma=\gamma_1$ and take $\gamma_l=0$ for $l \ge 2$ for the parameters
defined in Assumption $1^\prime$. Nevertheless, we note that there is a potential for the demonstration of further
improvements by Theorem \ref{theorem.asymp} and Corollary \ref{corollary.asym.linear} using $\gamma_l > 0$ for
$l=1$ and some $l \ge 2$ when $q$ is not a square (the situation is different when $q$ is a square, cf. \cite[Corollary 1]{T}).

For simplicity of notation, for
$\gamma=\gamma_1$ and $\gamma_l=0$ for $l \ge 2$, we denote the lower bounds of Theorem \ref{theorem.asymp} and
Corollary \ref{corollary.asym.linear} by $R_{\gamma,x_1, \ldots, x_m}(\delta)$ and $R_\gamma^{\rm lin}(\delta)$, respectively. In the examples below,
the required values of these two functions are computed by using Definition
\ref{definition.Psi} and Proposition
\ref{general.proposition.I.y.x1.x2.xm.compute}.

Let $R_{NO2,\gamma,x}(\delta)$ denote the lower bound in \cite[Theorem 5.1]{NO2}.
Moreover, let $R_{X,\gamma}^{\rm lin}(\delta)$ denote Xing's lower bound for $\alpha_q^{\rm lin}(\delta)$ in  \cite{X1} (see also \cite[Theorem 4.6]{NO2}).

\begin{examplerm} \label{example1}

Let $q=2^6$, $\gamma=\gamma_1=\sqrt{q}-1$, $\gamma_l=0$ for $l \ge 2$, and
\be
\begin{array}{l}
\dd \delta=\frac{13763868443250238929521503984833381597731412559044}{46065097831342932365531985486767649347321318605709} \\ \\
\dd =0.29879169026501515839 \ldots .
\end{array}
\nn\ee
In \cite[Example 5.2]{NO2}, using $x=10^{-13}$ it has been obtained that
\be
\alpha_q(\delta) \ge R_{NO2,\gamma,x}(\delta)=0.55835371587781529071 \ldots \; ,
\nn\ee
and it has been demonstrated that $R_{NO2,\gamma,x}(\delta) - R_{X,\gamma}^{\rm lin}(\delta) \ge 7.3387 \cdot 10^{-15}$.

By Corollary \ref{corollary.asym.linear} we obtain that
\be
\alpha_q^{\rm lin}(\delta) \ge R_\gamma^{\rm lin}(\delta)= 0.55835395724081743804 \ldots .
\nn\ee
Note that $R_\gamma^{\rm lin}(\delta)- R_{NO2,\gamma,x}(\delta) \ge 2.4136300214732 \cdot 10^{-7}$, and  $R_\gamma^{\rm lin}(\delta)$
is better than $R_{X,\gamma}^{\rm lin}(\delta)$. Hence we have an improvement on the lower bound for $\alpha^{\rm lin}_q(\delta)$
compared to  Xing's bound in \cite{X1}.

By Theorem \ref{theorem.asymp} with $x_1=3.41 \cdot 10^{-16}$, $x_2=1.0634 \cdot 10^{-23}$, and $x_3=1.93 \cdot 10^{-31}$,
we obtain
$\alpha_q(\delta) \ge R_{\gamma,x_1,x_2,x_3}(\delta)$, where
\be
R_{\gamma,x_1,x_2,x_3}(\delta) - R_\gamma^{\rm lin}(\delta) \ge 2.711029 \cdot 10^{-17}.
\nn\ee
Hence $R_{\gamma,x_1,x_2,x_3}(\delta)$ gives a further improvement on the lower bound for $\alpha_q(\delta)$.

Now let
\be
\begin{array}{l}
\dd \delta=\frac{32301229388092693436010481501934267749589906046665}{46065097831342932365531985486767649347321318605709} \\ \\
\dd =0.70120830973498484160 \ldots \; .
\end{array}
\nn\ee
In \cite[Example 5.2]{NO2}, using $x=10^{-13}$ it has been obtained that
\be
\alpha_q(\delta) \ge R_{NO2,\gamma,x}(\delta)=0.15593709640785805503 \ldots \; ,
\nn\ee
and it has been demonstrated that $R_{NO2,\gamma,x}(\delta) - R_{X,\gamma}^{\rm lin}(\delta) \ge 1.97862 \cdot 10^{-14}$.

By Corollary \ref{corollary.asym.linear} we obtain that
\be
\alpha_q^{\rm lin}(\delta) \ge R_\gamma^{\rm lin}(\delta)= 0.15593754394482448829 \ldots .
\nn\ee
Note that $R_\gamma^{\rm lin}(\delta)- R_{NO2,\gamma,x}(\delta) \ge 4.4753696643325 \cdot 10^{-7}$,  hence $R_\gamma^{\rm lin}(\delta)$
is better than $R_{X,\gamma}^{\rm lin}(\delta)$. Hence we have an improvement on the lower bound for $\alpha^{\rm lin}_q(\delta)$
compared to  Xing's bound in \cite{X1}.

By Theorem \ref{theorem.asymp} with $3.89\cdot 10^{-18}$, $x_2=1.98\cdot10^{-26}$, and $x_3=5.87\cdot10^{-35}$, we obtain
$\alpha_q(\delta) \ge R_{\gamma,x_1,x_2,x_3}(\delta)$, where
\be
R_{\gamma,x_1,x_2,x_3}(\delta) - R_\gamma^{\rm lin}(\delta) \ge 2.592642 \cdot 10^{-19}.
\nn\ee
Hence $R_{\gamma,x_1,x_2,x_3}(\delta)$ gives a further improvement on the lower bound for $\alpha_q(\delta)$.

\end{examplerm}

\begin{examplerm} \label{example2}

Let $q=7^2$, $\gamma=\gamma_1=\sqrt{q}-1$, $\gamma_l=0$ for $l \ge 2$, and
\be
\begin{array}{l}
\dd \delta=\frac{7334559589562321721169749749908497945081695123431}{18755194537338788993696079784908084949457099261873} \\ \\
\dd =0.39106816913897159912 \ldots \; .
\end{array}
\nn\ee
In \cite[Example 5.3]{NO2}, using $x=10^{-13}$ it has been obtained that
\be
\alpha_q(\delta) \ge R_{NO2,\gamma,x}(\delta)=0.44226734872224546020 \ldots \;,
\nn\ee
and it has been demonstrated that $R_{NO2,\gamma,x}(\delta) - R_{X,\gamma}^{\rm lin}(\delta) \ge 6.57561 \cdot 10^{-14}$.

By Corollary \ref{corollary.asym.linear} we obtain that
\be
\alpha_q^{\rm lin}(\delta) \ge R_\gamma^{\rm lin}(\delta)= 0.44226758374884970747 \ldots .
\nn\ee
Note that $R_\gamma^{\rm lin}(\delta)- R_{NO2,\gamma,x}(\delta) \ge 2.3502660424726 \cdot 10^{-7}$, and  $R_\gamma^{\rm lin}(\delta)$
is better than $R_{X,\gamma}^{\rm lin}(\delta)$. Hence we have an improvement on the lower bound for $\alpha^{\rm lin}_q(\delta)$
compared to  Xing's bound in \cite{X1}.

By Theorem \ref{theorem.asymp} with $x_1=1.93 \cdot 10^{-13}$, $x_2=1.53\cdot 10^{-19}$, and $x_3=7.08\cdot 10^{-26}$,
we obtain
$\alpha_q(\delta) \ge R_{\gamma,x_1,x_2,x_3}(\delta)$, where
\be
R_{\gamma,x_1,x_2,x_3}(\delta) - R_\gamma^{\rm lin}(\delta) \ge 1.857062 \cdot 10^{-14}.
\nn\ee
Hence $R_{\gamma,x_1,x_2,x_3}(\delta)$ gives a further improvement on the lower bound for $\alpha_q(\delta)$.

Now let
\be
\begin{array}{l}
\dd \delta=\frac{11420634947776467272526330034999587004375404138442}{18755194537338788993696079784908084949457099261873} \\ \\
\dd =0.60893183086102840087 \ldots \; .
\end{array}
\nn\ee
In \cite[Example 5.3]{NO2}, using $x=10^{-13}$ it has been obtained that
\be
\alpha_q(\delta) \ge R_{NO2,\gamma,x}(\delta)=0.22440368700019503856 \ldots \; ,
\nn\ee
and it has been demonstrated that $R_{NO2,\gamma,x}(\delta) - R_{X,\gamma}^{\rm lin}(\delta) \ge 7.21362 \cdot 10^{-14}$.

By Corollary \ref{corollary.asym.linear} we obtain that
\be
\alpha_q^{\rm lin}(\delta) \ge R_\gamma^{\rm lin}(\delta)= 0.22440401150099750683 \ldots .
\nn\ee
Note that $R_\gamma^{\rm lin}(\delta)- R_{NO2,\gamma,x}(\delta) \ge 3.2450080246826 \cdot 10^{-7}$, and  $R_\gamma^{\rm lin}(\delta)$
is better than $R_{X,\gamma}^{\rm lin}(\delta)$. Hence we have an improvement on the lower bound for $\alpha^{\rm lin}_q(\delta)$
compared to Xing's bound in \cite{X1}.

By Theorem \ref{theorem.asymp} with $x_1=5.86\cdot 10^{-14}$, $x_2=3.207\cdot 10^{-20}$, and $x_3=1.02\cdot 10^{-26}$, we obtain
$\alpha_q(\delta) \ge R_{\gamma,x_1,x_2,x_3}(\delta)$, where
\be
R_{\gamma,x_1,x_2,x_3}(\delta) - R_\gamma^{\rm lin}(\delta) \ge 5.258306 \cdot 10^{-15}.
\nn\ee
Hence $R_{\gamma,x_1,x_2,x_3}(\delta)$ gives a further improvement on the lower bound for $\alpha_q(\delta)$.

\end{examplerm}

\begin{examplerm} \label{example3}

Let $q=2^{21}$, $\dd \gamma=\gamma_1=\frac{2(q^{2/3}-1)}{q^{1/3}+2}$ (see (\ref{eq14})), $\gamma_l=0$ for $l \ge 2$, and
\be
\begin{array}{l}
\dd \delta=\frac{1034323484865452473463726110309814032498446010098}{99621193732964014413326435515634059733734238550355} \\ \\
\dd =0.01038256465424386359 \ldots \; .
\end{array}
\nn\ee
In \cite[Example 5.4]{NO2}, using $x=10^{-60}$ it has been obtained that
\be
\alpha_q(\delta) \ge R_{NO2,\gamma,x}(\delta)=0.98564990803085654665 \ldots \;,
\nn\ee
and it has been demonstrated that $R_{NO2,\gamma,x}(\delta) - R_{X,\gamma}^{\rm lin}(\delta) \ge 2.1335699248 \cdot 10^{-61}$.

By Corollary \ref{corollary.asym.linear} we obtain that
\be
\alpha_q^{\rm lin}(\delta) \ge R_\gamma^{\rm lin}(\delta)= 0.98564990803085654673 \ldots .
\nn\ee
Note that $R_\gamma^{\rm lin}(\delta)- R_{NO2,\gamma,x}(\delta) \ge 7 \cdot 10^{-20}$, and  $R_\gamma^{\rm lin}(\delta)$
is better than $R_{X,\gamma}^{\rm lin}(\delta)$. Hence we have an improvement on the lower bound for $\alpha^{\rm lin}_q(\delta)$
compared to  Xing's bound in \cite{X1}.

By Theorem \ref{theorem.asymp} with $x_1=6.29 \cdot10^{-65}$ and $x_2=7.09\cdot10^{-97}$,
we obtain
$\alpha_q(\delta) \ge R_{\gamma,x_1,x_2}(\delta)$, where
\be
R_{\gamma,x_1,x_2}(\delta) - R_\gamma^{\rm lin}(\delta) \ge 1.261672 \cdot 10^{-66}.
\nn\ee
Hence $R_{\gamma,x_1,x_2}(\delta)$ gives a further improvement on the lower bound for $\alpha_q(\delta)$.

Now let
\be
\begin{array}{l}
\dd \delta=\frac{98586870248098561939862709405324245701235792540257}{99621193732964014413326435515634059733734238550355} \\ \\
\dd =0.98961743534575613640 \ldots \; .
\end{array}
\nn\ee
In \cite[Example 5.4]{NO2}, using $x=10^{-60}$ it has been obtained that
\be
\alpha_q(\delta) \ge R_{NO2,\gamma,x}(\delta)=0.00641503733934427385 \ldots \; ,
\nn\ee
and it has been demonstrated that $R_{NO2,\gamma,x}(\delta) - R_{X,\gamma}^{\rm lin}(\delta) \ge 4.2225689802 \cdot 10^{-61}$.

By Corollary \ref{corollary.asym.linear} we obtain that
\be
\alpha_q^{\rm lin}(\delta) \ge R_\gamma^{\rm lin}(\delta)= 0.00641503733934427410 \ldots .
\nn\ee
Note that $R_\gamma^{\rm lin}(\delta)- R_{NO2,\gamma,x}(\delta) \ge 2.4 \cdot 10^{-19}$, and  $R_\gamma^{\rm lin}(\delta)$
is better than $R_{X,\gamma}^{\rm lin}(\delta)$. Hence we have an improvement on the lower bound for $\alpha^{\rm lin}_q(\delta)$
compared to  Xing's bound in \cite{X1}.

By Theorem \ref{theorem.asymp} with $x_1=6.5\cdot 10^{-86}$ and $x_2=2.4 \cdot 10^{-127}$,
we obtain
$\alpha_q(\delta) \ge R_{\gamma,x_1,x_2}(\delta)$, where
\be
R_{\gamma,x_1,x_2}(\delta) - R_\gamma^{\rm lin}(\delta) \ge 9.103449 \cdot 10^{-88}.
\nn\ee
Hence $R_{\gamma,x_1,x_2}(\delta)$ gives a further improvement on the lower bound for $\alpha_q(\delta)$.

\end{examplerm}

\end{document}